\documentclass[ejs]{imsart}

\RequirePackage[OT1]{fontenc}
\RequirePackage[numbers]{natbib}
\RequirePackage[colorlinks,citecolor=blue,urlcolor=blue]{hyperref}
\usepackage{amsmath}
\usepackage{amsfonts}
\usepackage{amsthm} 
\usepackage{amssymb}
\usepackage{amstext}
\usepackage{bbm}
\usepackage{mathtools}
\usepackage{graphicx}
\graphicspath{ {./figures/} }
\usepackage[mathscr]{euscript}
\usepackage{tabularx}
\usepackage{mathrsfs}
\usepackage{enumerate} 
\usepackage{cases} 
\usepackage{float} 
\usepackage{caption}
\usepackage{subcaption}
\usepackage{comment}
\usepackage{wrapfig}
\usepackage[utf8]{inputenc}
\usepackage[colorinlistoftodos]{todonotes}
\usepackage{amsbsy} 
\usepackage{physics} 
\usepackage{upgreek} 
\usepackage{diagbox} 
\usepackage{multirow}
\usepackage{booktabs}
\usepackage[normalem]{ulem}

\theoremstyle{plain}

\newtheorem{thm}{Theorem}[section]
\newtheorem{lem}[thm]{Lemma}

\theoremstyle{definition}

\theoremstyle{remark}

\def\bbE{\mathbbm {E}}

\def \Dual {\mathcal {D}}

\def \bsm {\boldsymbol}
\def \btheta {\boldsymbol{\theta}} 

\def \bpsi {\boldsymbol{\psi}}
 
\def \blambda {\boldsymbol{\lambda}} 
 
\def \bnu {\boldsymbol{\nu}}

\def \bxi {\boldsymbol{\xi}}

\def \bq {\boldsymbol{\mathrm{q}}} 
\def \bh {\boldsymbol{\mathrm{h}}}

\def \bone {\boldsymbol{\mathrm{1}}} 
\def \be {\boldsymbol{e}}

\def \ba {\boldsymbol{a}}

\def \by {\boldsymbol{y}}

\newcommand{\ind}{\mbox{$\mathbbm{1}$}}

\begin{document}

\begin{frontmatter}

\title{Empirical likelihood ratio test on quantiles under a density ratio model}
\runtitle{ELRT on quantiles under a DRM}

\begin{aug}
\author[UBC]{\fnms{Archer Gong} \snm{Zhang}\ead[label=e1]{archer.zhang@stat.ubc.ca}},
\author[URI]{\fnms{Guangyu} \snm{Zhu}\ead[label=e2]{guangyuzhu@uri.edu}}
\and 
\author[UBC]{\fnms{Jiahua} \snm{Chen}\ead[label=e3]{jhchen@stat.ubc.ca}}

\address[UBC]{Department of Statistics \\
The University of British Columbia \\ 
Vancouver, BC, Canada, V6T 1Z4 \\ 
\printead{e1,e3}} 

\address[URI]{Department of Computer Science and Statistics \\
The University of Rhode Island \\ 
Kingston, RI, USA 02881 \\
\printead{e2}} 

\runauthor{A. Zhang, G. Zhu, and J. Chen}
\end{aug} 

\begin{abstract}
Population quantiles are important parameters in many applications. 
Enthusiasm for the development of effective statistical inference procedures for quantiles and their functions has been high for the past decade. 
In this article, we study inference methods for quantiles when multiple samples from linked populations are available. 
The research problems we consider have a wide range of applications. 
For example, to study the evolution of the economic status of a country, economists monitor changes in the quantiles of annual household incomes, based on multiple survey datasets collected annually. 
Even with multiple samples, a routine approach would estimate the quantiles of different populations separately. 
Such approaches ignore the fact that these populations are linked and share some intrinsic latent structure. 
Recently, many researchers have advocated the use of the density ratio model (DRM) to account for this latent structure and have developed more efficient procedures based on pooled data. 
The nonparametric empirical likelihood (EL) is subsequently employed. 
Interestingly, there has been no discussion in this context of the EL-based likelihood ratio test (ELRT) for population quantiles. 
We explore the use of the ELRT for hypotheses concerning quantiles and confidence regions under the DRM. 
We show that the ELRT statistic has a chi-square limiting distribution under the null hypothesis. 
Simulation experiments show that the chi-square distributions approximate the finite-sample distributions well and lead to accurate tests and confidence regions. 
The DRM helps to improve statistical efficiency.
We also give a real-data example to illustrate the efficiency of the proposed method. 
\end{abstract}



%\tableofcontents 

\end{frontmatter}

%Main body of the paper 
\section{Introduction} 
\label {sec:intro}

Suppose we have $ m+1 $ independent random samples from population distributions $ G_{0}, G_{1}, \ldots, G_{m} $. 
Let their respective density functions with respect to some $ \sigma $-finite measure be $ g_{k} (\cdot) $. 
If there exist a vector-valued function $ \bq (x) $ and unknown vector-valued parameters $ \btheta_{k} $ of dimension $ d $ such that
\begin{align} 
\label {DRM}  
g_{k} (x) = \exp \{ \btheta_{k}^{\top} \bq (x) \} g_{0} (x), 
\end{align} 
then they define a density ratio model (DRM) as introduced by \citet {anderson1979multivariate}. 
By convention, we call $ G_{0} $ the base distribution and $ \bq (x) $ the basis function. 
There is a symmetry in the DRM: any one of $ G_{0}, \ldots, G_{m} $ may serve as the base distribution. 
We require the first element of $ \bq (x) $ to be 1 so that the corresponding coefficient is a normalization constant, and the elements of $ \bq (x) $ must be linearly independent. 
The linear independence is a natural requirement: otherwise, some elements of $ \bq (x) $ are redundant. 

When data are collected from a DRM, the whole data set can be utilized to estimate 
$ G_{0} $, which lead to efficiency gain.
The nonparametric $ G_{0} $ assumption in the DRM is nonrestrictive. 
Combined with a moderate-sized $ \bq (x) $, a single DRM contains a broad range of parametric distribution families. 
Thus, the DRM has a low risk of model misspecification. 
There is a growing interest in the DRM in statistics \citep {qin1998inferences, fokianos2001semiparametric, de2017bayesian, zhuang2019semiparametric} as well as in the machine learning community \citep {sugiyama2012density}.
In this paper, we study the inference problem for population quantiles under the DRM. 
Population quantiles and their functions are important parameters in many applications.
For example, government agents gauge the overall economic status of a country based on annual surveys of household income distribution.
The trend in the quantiles of the income distribution is indicative \citep {berger2003variance, muller2008measurement}.
In forestry, the lower quantiles of the mechanical strength of wood products are vital design values \citep {verrill2015simulation}. 
Other examples include 
\citet{chen1993smoothed, yang2012bayesian, chen2013quantile, chen2016monitoring, koenker2017handbook, gonccalves2018dynamic, chen2019small}. 

The data from DRMs are a special type of biased sample 
\citep{vardi1982nonparametric, vardi1985empirical, qin1998inferences, qin2017connections}.
The empirical likelihood (EL) of \citet {owen2001empirical} is an ideal platform for statistical inference under the DRM. 
The EL retains the effectiveness of likelihood methods and does not impose a restrictive parametric assumption. 
The ELRT statistic has a neat chi-square limiting distribution, 
much like the parametric likelihood ratio test given independent 
and identically distributed (i.i.d.) observations 
\citep {owen1988empirical, qin1994empirical}.
The EL has already been widely used for data analysis under the DRM 
\citep {qin1993empirical, qin1997goodness,chen2013quantile, cai2017hypothesis}.
However, there has been limited discussion of the ELRT in the biased sampling context.
Both \citet {qin1993empirical} and \citet {cai2017hypothesis} permit no additional equations.
Although the classical Wald method remains effective for both hypothesis tests 
and confidence regions \citep{qin1998inferences, chen2013quantile, chen2016monitoring}, 
it must be aided by a consistent and stable variance estimate. 
In addition, its confidence regions are oval-shaped regardless of the shape of the data cloud.
Thus, an ELRT has the potential to push the boundary of the DRM much further.

This paper establishes the limiting chi-square distribution of the ELRT for quantiles under the DRM. 
We prove that the ELRT statistic has a chi-square limiting distribution under certain conditions. 
The resulting confidence regions have data-driven shapes, more accurate coverage probabilities, and smaller volumes. 
In Section \ref {sec:ELRT}, we state the problem of interest and the proposed ELRT under the DRM. 
In Section \ref {sec:main_results}, we study the limiting distribution of the ELRT statistic and some other useful asymptotic results. 
We illustrate the superiority of the ELRT and the associated confidence regions through simulated data in Section \ref {sec:simulations} and for real-world data in Section \ref {sec:realdata}. 
We illustrate the power property of the ELRT in Section \ref{power}.
Technical details and the proofs of the main theorems are given in Appendices \ref {App:Proofs} and \ref {App:define_profile}.

%%%%%%%%%%%%%%%%%%%%%%%%%%%%%%%%%%%%%%%%%%%%%%%%%%%%%%%%%%%%%%%%%%%%%%

\section{Research problem and proposed approach} 
\label {sec:ELRT}

Let $ \{ x_{k j}: 1 \leq j \leq n_{k}, 0 \leq k \leq m  \} $ be $ m+1 $ independent i.i.d. samples from a DRM defined by \eqref {DRM}. 
Let $ n = \sum_{k=0}^m n_k $ be the total sample size. 
Denote by $ \xi_{k} $ the $ \tau_{k} $ quantile of the $ k $th population for some $ \tau_{k} \in (0, 1) $ and $ k = 0, 1, \ldots, m $. 
Let $ \bxi = \{ \xi_{k}: k \in I \} $ be the quantiles at some levels of populations in an index set $ I \subseteq \{ 0, 1, \ldots, m \} $ of size $ l $. 
We study the ELRT under the DRM for the following hypothesis:
\begin{align} 
H_{0}: \bxi = \bxi^{*} \, \, \, \text { against } \, \, \, H_{1}: \bxi \neq \bxi^{*}, 
\label {simple_hyp}   
\end{align}  
for some given $ \bxi^{*} $ of dimension $ l $. 

The hypothesis formulated in \eqref {simple_hyp} has many applications. 
In socio-economic studies, 
when studying the distributions of household disposable incomes, economists 
and social scientists often divide the collected survey data into five groups. 
These groups are famously known as quintile groups. 
The first group consists of the lowest $ 20\% $ of the data,
the second group consists of the next 20\%, 
and so on. 
Many studies have shown that the quintiles are important for 
explaining the economy and consumer behaviour 
\citep {castello2002human, wunder2012income, humphries2014households, corak2018canadian}. 
In statistics, the cut-off points of these quintile groups are the quantiles of the populations: 
for example, the $ 20 $th percentile separates the first and second quintile groups. 
Governments may, therefore, consider this $ 20 $th percentile as key for 
determining which families should receive a special subsidy to help society's less fortunate. 
Moreover, when new policies are implemented, the evolution of the 
quantiles of household income over time may reflect the impact of the policies. 
As a consequence, these quantiles are of particular interest to 
social scientists and politicians as a way to measure the effects of policy changes.
In statistical inference, these types of tasks can most appropriately be carried out 
using a hypothesis testing procedure, 
which can be naturally extended to the construction of confidence regions. 
Hence, the research problem we study here is of scientific significance in many applications. 
In the real-data analysis, we study confidence regions for quantiles 
of household incomes based on US Consumer Expenditure Surveys. 

We use an ELRT to test the hypothesis in \eqref {simple_hyp}. 
Let $ p_{k j} = \mathrm {d} G_{0} (x_{k j}) = P (X = x_{k j}; G_{0}) $ for all applicable $ k, j $. 
The EL function is the probability of observing the data. 
Under the DRM, it is given by 
\begin{align} 
L_{n} (G_{0}, \ldots, G_{m}) 
= 
\prod_{k, j} \mathrm {d} G_{k} (x_{k j}) 
 = \big \{ \prod_{k, j} p_{k j} \big \} \big \{ \prod_{k, j} \exp (\btheta_{k}^{\top} \bq (x_{k j})) \big \}. 
\label {EL}
\end{align}
For notational convenience, we have dropped the ranges of the indices in the expressions. 
Observe that the EL in \eqref {EL} is 0 if $ G_{0} $ is a continuous distribution. 
Surprisingly, this seemingly devastating property does little harm to the usefulness of the EL.
Since the EL in \eqref {EL} can also be regarded as a function of the parameters 
$ \btheta \coloneqq \{ \btheta_{r}: 1 \leq r \leq m \} $ 
and the base distribution $ G_{0} $, we may write its logarithm as 
\begin{align*} 
\ell_{n} ( \btheta, G_{0} ) 
= 
\log L_{n} (G_{0}, \ldots, G_{m}) = \sum_{k, j} \log p_{k j} 
+ 
\sum_{k, j} \btheta_{k}^{\top} \bq (x_{k j}), 
\end{align*} 
where we define $ \btheta_{0} = \bsm {0} $ by convention. 

Let $ \mathbbm {E}_{r} $ be the expectation operation under $ G_{r} $,
and let
\[
h_{r} (x, \btheta) 
= 
\exp (\btheta_{r}^{\top} \bq (x))
\]
be the density of $ G_{r} $ with respect to $ G_{0} $ for $ r = 0, 1, \ldots, m $. 
Clearly, $ h_{0} (x, \btheta) = 1 $.
This also implies that
\begin{align} 
\label{auto.eq}
\mathbbm {E}_{0} [ h_{r} (X, \btheta) ] = \mathbbm {E}_{0} \left [ \exp (\btheta_{r}^{\top} \bq (X)) \right ] = 1.
\end{align} 
The $ \tau_{r} $ population quantile $ \xi_{r} $ of $ G_{r} $ satisfies or is defined to be a solution of
\begin{align} 
\label{quantile.eq}
\mathbbm {E}_{r} \big [ \mathbbm{1} (X \leq \xi_{r}) - \tau_{r} \big ] 
= 
\mathbbm {E}_{0} 
\left [ h_{r} (X, \btheta) \{ \mathbbm {1} (X \leq \xi_{r}) - \tau_{r} \} \right ] = 0.
\end{align} 
Let
\[
\varphi_{r} (x, \btheta, \bxi) 
= 
h_{r} (x, \btheta) [ \mathbbm {1} (x \leq \xi_{r}) - \tau_{r} ].
\]
Following \citet {owen2001empirical} and \citet {qin1994empirical}, we introduce the profile log-EL of the population quantiles $ \bxi $:
\begin{align} 
\label{profile.eq}
\tilde \ell_{n} (\bxi) 
= 
\sup_{\btheta, G_{0}} \Big \{ \ell_{n} (\btheta, G_{0}) ~ | ~
& \sum_{k, j} p_{k j} h_{r} (x_{k j}, \btheta)  = 1, \, r = 0, 1, \ldots, m,
\nonumber \\ 
& \sum_{k, j} p_{k j} \varphi_{r} (x_{k j}, \btheta, \bxi) = 0, \, r \in I \Big \}
\end{align} 
and
\[
\sup_{ \btheta, G_{0} } \{ \ell_{n} (\btheta, G_{0}) \} 
=
\sup_{ \btheta, G_{0} } \{ \ell_{n} (\btheta, G_{0}) | 
\sum_{k, j} p_{k j} h_{r} (x_{k j}, \btheta) = 1, r = 0, 1, \ldots, m \}.
\]
An ELRT statistic for the hypothesis in \eqref {simple_hyp} is defined as
\begin{align*} 
R_{n} = 2 \left [ \sup_{ \btheta, G_{0} } \{ \ell_{n} (\btheta, G_{0})\} - \tilde \ell_{n} (\bxi^{*}) \right ]. 
\end{align*}

We call $ R_{n} $ the ELRT statistic hereafter. 
Clearly, the larger the value of $ R_{n} $, the stronger the evidence for departure from the null hypothesis in the direction of the alternative hypothesis. 
We reject $ H_{0} $ when $ R_{n} $ exceeds some critical value that is decided based on the distributional information of $ R_{n} $ under $ H_{0} $. 
The limiting distribution of $ R_{n} $ and other related properties are given in the next section.

We observe that the approach needs no change for a set of quantiles from the same population. 
For notational simplicity, the presentation is given for quantiles from different populations.

%%%%%%%%%%%%%%%%%%%%%%%%%%%%%%%%%%%%%%%%%%%%%%%%%%%%%%%%%%%%%%%%%%%%%%

\section{Asymptotic properties of $ R_{n} $ and other quantities} 
\label{sec:main_results}

The distributional information of $ R_{n} $ is vital to the implementation of the ELRT in applications. 
In this section, we show that it is asymptotically chi-square distributed. 
We also present some secondary but useful asymptotic results.

\subsection{A dual function}
 
The profile log-EL function $ \tilde \ell_{n} (\bxi^{*}) $ is defined to be the solution of an optimization problem that can be solved by the Lagrange multiplier method. 
Let $ \boldsymbol {t} = (t_{0}, \ldots, t_{m})$ and $\blambda = \{ \lambda_{r}: r \in I \} $ be Lagrange multipliers. 
Define a Lagrangian as 
\begin{align*} 
\mathcal{L} (\boldsymbol{t}, \blambda, \btheta, G_{0}) 
= & \ell_{n} (\btheta, G_{0})
+ \sum_{r = 0}^{m} t_{r} \big \{ 1 - \sum_{k, j} p_{k j} h_{r} (x_{k j}, \btheta) \big \}
\nonumber \\ 
& 
- \sum_{r \in I} n \lambda_{r} \big \{ \sum_{k, j} p_{k j} \varphi_{r} (x_{k j}, \btheta, \bxi^{*}) \big \}.   
\end{align*} 

In Appendix \ref {App:define_profile}, we will show that under mild conditions that are easy to verify, there aways exists some $ \btheta $ such that a solution in $ G_{0} $ to \eqref {auto.eq} and \eqref {quantile.eq} exists.
With this promise, according to the Karush--Kuhn--Tucker theorem \citep {boyd2004convex}, the solution to the constrained optimization problem in \eqref{profile.eq} satisfies 
\begin{align*} 
\frac { \partial \mathcal{L} (\boldsymbol{t}, \blambda, \btheta, G_{0}) } 
       { \partial (\boldsymbol {t}, \blambda, \btheta, p_{k j}) } 
       = \bsm {0}.  
\end{align*} 
Let
$ (\hat { \boldsymbol {t} }, \hat \blambda, \hat \btheta, \hat p_{k j}) $ be the solution. 
Some simple algebra gives $ \hat t_{r} = n_{r} $ and 
\begin{align*} 
\hat p_{k j}
= n^{-1} \left \{ 
\sum_{r = 0}^{m} \rho_{r} h_{r} (x_{k j}, \hat \btheta) 
+ 
\sum_{r \in I} \hat \lambda_{r} \varphi_{r} (x_{k j}, \hat \btheta, \bxi^{*}) \right \}^{-1}, 
\end{align*} 
where $ \rho_r = n_{r}/n $.  

We now introduce another set of notation: 
\begin{gather*} 
\bar {h} (x, \btheta) 
 = \sum_{r = 0}^{m} \rho_{r} h_{r} (x, \btheta), \\
\bh (x, \btheta) 
 = 
(\rho_{1} h_{1} (x, \btheta)/\bar {h} (x, \btheta), \ldots, \rho_{m} h_{m} (x, \btheta)/\bar {h} (x, \btheta))^{\top}, \\
\psi_{r} (x, \btheta) = \varphi_{r} (x, \btheta, \bxi^{*})/\bar {h} (x, \btheta), \\ 
\bpsi (x, \btheta) = \{ \psi_{r} (x, \btheta) : r \in I \}.
\end{gather*} 
To aid our memory, we note that
$ \bar {h} (x, \btheta) $ is a mixture density with mixing proportions $ \rho_{0}, \ldots, \rho_{m} $; 
$ \bh (x, \btheta) $ is a vector of density functions with respect to the mixture $ \bar {h} (x, \btheta) $ combined with the mixing proportions; 
and $ \bpsi (x, \btheta) $ is a vector of normalized $ \varphi_{r} (x, \btheta, \bxi^{*}) $.
With the help of this notation, we define a {\em dual function}
\begin{align}
\Dual (\blambda, \btheta) 
= &
\sum_{k, j} \btheta_{k}^{\top} \bq (x_{k j})
 -
\sum_{k, j} \log \bar {h} (x_{k j}, \btheta) 
\nonumber \\
& - \sum_{k, j} 
\log \big \{ 1 + \sum_{r \in I} \lambda_{r} \psi_{r} (x_{k j}, \btheta) \big \}.
\label {new_dual1} 
\end{align}
The dual function has some easily verified mathematical properties.
We can show that
\begin{align} 
\tilde \ell_{n} (\bxi^{*}) = \Dual (\hat \blambda, \hat \btheta) - n \log n,
\label {dual_relation}
\end{align} 
and that $ (\hat \blambda, \hat \btheta) $ is a saddle point of $ \Dual (\blambda, \btheta) $ satisfying
\begin{align} 
\frac {\partial \Dual (\blambda, \btheta)} {\partial (\blambda, \btheta)} = \bsm {0}. 
\label {joint_solution}
\end{align}  

In the following section, we study some of the properties of $ \tilde \ell_{n} (\bxi^{*}) $ through the dual function $ \Dual (\blambda, \btheta) $.

\subsection{Asymptotic properties} 

We discuss the asymptotic properties under the following nonrestrictive conditions on the sampling plan and the DRM.

\noindent
{\bf Conditions:} 
\begin{enumerate}[(i)] 
\item 
\label{Condition.i}
The sample proportions $ \rho_{k} = n_{k}/n $ have limits in $ (0, 1) $ as $ n \to \infty $;

\item 
\label{Condition.ii}
The matrix $ \bbE_{0} [ \bq (X) \bq^{\top} (X) ] $ is positive definite; 

\item 
\label{Condition.iii}
For each $ k = 0, 1, \ldots, m $ and $ \btheta_{k} $ in a neighbourhood of the true parameter value $ \btheta_{k}^{*} $, we have
\[
\bbE_{0} \left [ \exp (\btheta_{k}^{\top} \bq (X)) \right ]
 =
\bbE_{0} [ h_{k} (X, \btheta) ] < \infty. 
\]
\end{enumerate}  

Here are some implications of the above conditions.

\begin{enumerate} 

\item 
Under Condition \eqref{Condition.iii}, the moment generating function of $ \bq (X) $ with respect to $ G_{k} $ exists in a neighbourhood of $ \bsm {0} $.
Hence, all finite-order moments of $ \| \bq (X) \| $ are finite. 

\item 
When $ n $ is large enough and $ (\blambda, \btheta) $ is in a small neighbourhood of $ (\bsm {0}, \btheta^{*}) $, the derivatives of the dual function $ \Dual (\blambda, \btheta) $ are all bounded by some polynomials of $ \| \bq (x) \|$.
Hence, they are all integrable. 

\item 
Under Condition \eqref{Condition.ii}, the sample version of $ \bbE_{0} [\bq (X) \bq^{\top} (X)] $ is also positive definite when $ n $ is very large. 

\end{enumerate} 

We now state the main results; 
the proofs are given in Appendix \ref {App:Proofs}. 

\begin{lem} 
\label{lem:second_deriv}

Under Conditions \eqref{Condition.i} to \eqref{Condition.iii}, as $ n \to \infty $, 
\[ 
n^{-1} \left. \frac {\partial^{2} \Dual (\blambda, \btheta)} {\partial (\blambda, \btheta) \partial (\blambda, \btheta)^{\top}} 
\right \rvert_{\blambda = \bsm {0}, \btheta = \btheta^*} 
\to 
S, 
\] 
almost surely for some full-rank square matrix $ S $ of dimension $ (dm+l) $ that has the expression
\[ 
S =
\sum_{k = 0}^{m} \rho_{k} \mathbbm {E}_{k} \left [ \frac { \partial^{2} \Dual_{k} (X, \bsm {0}, \btheta^{*}) } { \partial (\blambda, \btheta) \partial (\blambda, \btheta)^{\top} } \right ].
\]

\end{lem} 

The second derivative of the dual function $ \Dual (\blambda, \btheta) $ is not negative definite in comparison to a usual likelihood function. 
This is understandable because $ \blambda $ is not a model parameter.
However, it has full rank and plays an important role in localizing $ \hat \btheta $.  

The next result implies that the dual function $ \Dual (\blambda, \btheta) $ resembles the log-likelihood function under regularity conditions in an important way: its first derivative is an unbiased estimating function.

\begin{lem} 
\label{lem:first_deriv}
Under Conditions \eqref{Condition.i} to \eqref{Condition.iii}, we have
\[ 
\mathbbm {E} 
\left [
\left. \frac{\partial \Dual (\blambda, \btheta) }{\partial (\blambda, \btheta)} \right ] 
\right \rvert_{\blambda = \bsm {0}, \btheta = \btheta^*} = \bsm {0},
\] 
where the expectation is calculated by regarding $ x_{k j} $ as a random variable with distribution $ G_{k} $.

Furthermore, as $ n \to \infty $, we have
\[
 n^{-1/2}
\left. \dfrac {\partial \Dual (\blambda, \btheta)} {\partial (\blambda, \btheta)} \right \rvert_{\blambda = \bsm {0}, \btheta = \btheta^*}
\overset {d} {\to}
N (\bsm {0}, V),  
\]
where $ V $ is a square matrix of dimension $ (dm+l) $.  

\end{lem}

A key step in the asymptotic study of $ \hat \btheta $ and the ELRT statistic $ R_{n} $ is localization. 
That is, $ \hat \btheta $ is in a small neighbourhood of the true value $ \btheta^{*} $ as the sample size $ n $ goes to infinity. 
The following lemma asserts that $ \hat \btheta $ is almost surely located in the $ O (n^{-1/3}) $-neighbourhood of $ \btheta^{*} $.  

\begin{lem} 
\label{lem:para_normality} 

Under Conditions \eqref{Condition.i} to \eqref{Condition.iii}, as $ n \to \infty $, the saddle point $ (\hat \blambda, \hat \btheta) $ of the dual function $ \Dual (\blambda, \btheta) $ is in the $ n^{-1/3} $-neighbourhood of $ (\bsm {0}, \btheta^{*}) $ with probability $ 1 $.

In addition, $ \sqrt {n} (\hat \blambda, \hat \btheta - \btheta^{*}) $ is asymptotically multivariate normal. 

\end{lem}

The results in the previous lemma shed light on the asymptotic properties of the EL under the DRM. 
At the same time, they pave the way for the following celebrated conclusion in the EL literature.

\begin{thm} 
\label{thm:ELRT_Stat_chisq}

Under Conditions \eqref {Condition.i} to \eqref {Condition.iii} and the null hypothesis \eqref{simple_hyp}, as $ n \to \infty $, the ELRT statistic 
\[
R_{n}  = 2 \left [ \sup_{ \btheta, G_{0} } \{ \ell_{n} (\btheta, G_{0}) \} - \tilde \ell_{n} (\bxi^{*}) \right ] 
\overset {d} \to \chi_{l}^{2}.
\]

\end{thm} 

Theorem \ref {thm:ELRT_Stat_chisq} enables us to determine an approximate rejection region for the ELRT. 
We reject the null hypothesis at the significance level $ \alpha $ when the observed value of $ R_{n} $ is larger than the upper $ \alpha $ quantile of the chi-square distribution $ \chi_{l}^{2} $.
This also provides a foundation for the construction of confidence regions of $ \bxi $.
Let
\[ 
R_{n} (\bxi) 
= 2 \left [ \sup_{ \btheta, G_{0} } \{ \ell_{n} (\btheta, G_{0}) \} - \tilde \ell_{n} (\bxi) \right ]. 
\] 
An ELRT-based $ (1 - \alpha) $ approximate confidence region for $ \bxi $ is
\begin{align} 
\{ \bxi: R_{n} (\bxi) \leq \chi_{l}^{2} (1 - \alpha) \}, 
\label {ELRT_CR}
\end{align} 
where $ \chi^{2}_{l} (1 - \alpha) $ is the $ (1 - \alpha) $ quantile of $ \chi_{l}^{2} $.

%%%%%%%%%%%%%%%%%%%%%%%%%%%%%%%%%%%%%%%%%%%%%%%%%%%%%%%%%%%%%%%%%%%%%%

\section{Simulation studies} 
\label{sec:simulations}

In this section, we report some simulation results.
We conclude that the chi-square approximation to the sample distribution of $ R_{n} $ is very accurate.
The corresponding confidence regions have a data-driven shape and accurate coverage probabilities.
In almost all cases considered, the $ R_{n} $-based confidence regions outperform 
those based on the Wald method in terms of the average areas and coverage probabilities.
The DRM markedly improves the statistical efficiency, and
the details are as follows.

\subsection{Numerical implementation and methods included} 

Recall that the ELRT statistic $ R_{n} $ is defined to be 
\[ 
R_{n} = 2 \left [ \sup_{\btheta, G_{0}} \{ \ell_{n} (\btheta, G_{0}) \} - \tilde \ell_{n} (\bxi^{*}) \right ]. 
\] 
In data analysis, we must solve the optimization problem 
$ \sup_{\btheta, G_{0}} \{ \ell_{n} (\btheta, G_{0}) \} $. 
As \citet {cai2017hypothesis} suggest, 
it can be transformed into an optimization problem of a convex function, 
and it has a simple solution. 
We further turn this optimization problem into the problem of solving 
a system of equations that are formed by equating the derivatives 
of the induced convex function to $ \bsm {0} $. 
The numerical implementation can be efficiently carried out by 
a root solver in the {\tt R} \citep {R} package {\tt nleqslv} \citep {nleqslv} for 
nonlinear equations.
It uses either the Newton or Broyden iterative algorithms.

To compute $ \tilde \ell_{n} (\bxi^{*}) $, we can solve \eqref {joint_solution}, 
as \eqref {dual_relation} suggests.
This leads to a system of $ dm+l $ nonlinear equations in $ (\blambda, \btheta) $, 
with $ d $ being the dimension of the vector-valued basis function $ \bq (x) $ 
and $ l $ the number of population quantiles of interest specified in $ \bxi^{*} $. 
In most applications, a $ \bq (x) $ with dimension 4 or less is suitable.
For a system of this size, the {\tt R} package {\tt nleqslv} for roots is very effective 
even when $ m $ is as large as $ 20 $. 
The existence of the solution to \eqref {auto.eq} and \eqref {quantile.eq} is proved in 
Appendix \ref {App:define_profile}.
Guided by this proof, our choice of the initial $ \blambda $ and $ \btheta $ 
guarantees numerical success. 

As is typical for DRM examples, we simulate data from the normal and 
gamma distributions and examine the ELRT-based hypothesis tests and 
confidence regions for the population quantiles. 
For comparison, we include Wald-based and nonparametric inference on the same quantiles. 
To make the article self-contained, we now briefly review the Wald 
and nonparametric methods.    
~ \\ 

\noindent 
{\bf Wald method.}  
The Wald method for confidence region construction of $ \bxi $ was given in \citet {chen2013quantile}. 
Let $ (\tilde \btheta, \tilde G_{0}) $ be the argument maximizer of $ \sup_{ \btheta, G_{0} } \{ \ell_{n} (\btheta, G_{0}) \} $, and also let 
\[ 
\tilde G_{r} (x) 
= \sum_{k, j} \mathbbm {1} (x_{k j} \leq x) h_{r} (x_{k j}, \tilde \btheta) \mathrm {d} \tilde G_{0} (x_{k j}), 
\] 
for $ r = 1, \ldots, m $, where $ \mathrm {d} \tilde G_{0} (x) = \tilde G_{0} (x) - \tilde G_{0} (x_{-}) $. 
The maximum EL estimator (MELE) of the $ \tau_{r} $ quantile of $ G_{r} $ is then given by 
\[ 
\tilde \xi_{r} = \inf \{ x: \tilde G_{r} (x) \geq \tau_{r} \}. 
\] 
Let $ \tilde \bxi = \{ \tilde \xi_{r}: r \in I \} $.
We have, as $ n \to \infty $, 
\[ 
\sqrt {n} (\tilde \bxi - \bxi^{*}) \to N (\bsm {0}, \Omega), 
\] 
for some matrix $ \Omega $ that is a function of $ G_{r} $ and $ \btheta $. 
A plug-in estimate $ \tilde \Omega $ of $ \Omega $ was suggested by \citet {chen2013quantile}, 
and an {\tt R} package {\tt drmdel} \citep {drmdel} by the authors of \citet {cai2017hypothesis} includes the MELE $ \tilde \bxi $ and $ \tilde \Omega $ in its output.
A level $ (1 - \alpha) $ approximate confidence region for $ \bxi $ based on the Wald method is then given by 
\begin{align} 
\{ \bxi: n (\tilde \bxi - \bxi)^{\top} \tilde \Omega^{-1} (\tilde \bxi - \bxi)
\leq 
\chi_{l}^{2} (1 - \alpha) \}. 
\label {Wald_CR}
\end{align} 
The Wald method can also be used for hypothesis tests on quantiles.
We refer to the confidence region in \eqref {Wald_CR} as the one based on the {\em Wald method}.

~ \\ 
\noindent 
{\bf Nonparametric method.}  
Suppose $ \hat G_{r} (x) = n_{r}^{-1} \sum_{j = 1}^{n_{r}} \ind (x_{r j} \leq x) $ 
is the empirical distribution based on a sample from the distribution $ G_{r} $, 
and $ \hat \xi_{r} $ is the sample quantile. 
The sample quantile is asymptotically normal \citep {serfling2000approximation} 
with asymptotic variance $ {\tau_{r} (1 - \tau_{r})}/{(\rho_{r} g_{r}^{2} (\xi_{r}))} $ 
as $ n \to \infty $ and $ n_{r}/n \to \rho_{r} $. 
In view of this, the Wald method remains applicable with the help of a 
nonparametric consistent density estimator.
We follow the literature and let
\[ 
\hat g_{r} (x) = \frac {1} {n_{r} b_{r}} \sum_{j = 1}^{n_{r}} K \left ( \frac {x_{r j} - x} {b_{r}} \right ), 
\] 
for some kernel function $ K (\cdot) $ and bandwidth $ b_{r} $. 
Under mild conditions on $ g_{r} (\cdot) $ and proper choices of $ K (\cdot) $ 
and $ b_{r} $, $ \hat g_{r} (x) $ is consistent \citep {silverman1986density}. 
We set $ K (\cdot) $ to the density function of the standard normal distribution, 
and we use a rule-of-thumb bandwidth suggested by \citet {silverman1986density}: 
\[ 
b_{r} = 0.9 \min \{ \hat \sigma_{r}, \widehat {\mathrm {IQR}}_{r}/1.34 \} n_{r}^{-1/5},  
\]
where $ \hat \sigma_{r}$ is the standard deviation of $ \hat G_{r} $ 
and $\widehat {\mathrm {IQR}}_{r}$ is the interquartile range. 
With these, we obtain a plug-in estimate
\[
\hat T \coloneqq \mathrm {diag} \{ \tau_{r} (1 - \tau_{r})/(\rho_{r} \hat g_{r}^{2} (\hat \xi_{r})): r \in I \}, 
\]
and subsequently a $ (1 - \alpha) $ approximate confidence region for $ \bxi $:
\begin{align} 
\{ \bxi: 
n (\hat \bxi - \bxi)^{\top} \hat T^{-1} (\hat \bxi - \bxi)
\leq 
\chi_{l}^{2} (1 - \alpha) \},  
\label {Nonpara_CR}
\end{align} 
where $ \hat \bxi = \{ \hat \xi_{r}: r \in I \} $.
This nonparametric Wald method can also be employed for hypothesis tests on quantiles.
We refer to the confidence region in \eqref {Nonpara_CR} 
as the one based on the {\em nonparametric method}.  

When constructing the confidence region in \eqref {Nonpara_CR}, 
density estimation is required as an intermediate step to obtain a variance estimate $ \hat T $. 
One may also use bootstrap method as an alternative nonparametric method 
to construct confidence regions for quantiles. 
We do not think these two nonparametric methods will lead to significantly different results, 
and hence we use \eqref {Nonpara_CR} as a nonparametric competitor in this article.

The proposed ELRT method apparently has the highest computational cost,
yet it takes a negligible second for each simulation repetition.
This renders recording the computational times unnecessary in simulation.

\subsection{Data generated from normal distributions} 

Normality is routinely assumed but unlikely strictly valid in real-world applications.
When multiple samples are available, we include all normal distributions without 
a normality assumption via a DRM coupled with $ \bq (x) = (1, x, x^{2})^{\top} $.
In this simulation, we generate data from $ m+1 = 6 $ normal distributions with sample sizes $ n_{r} = 100 $. 
Their means and standard deviations are chosen to be $ (0, 0, 1, 1, 2, 2) $ and $ (1, 1.2, 1.3, 1.5, 2, 1.5) $.
In the simulation experiment, we generate 1000 sets of samples of size $ n_{r} = 100 $ 
and compute the $ R_{n} $ values for the hypothesis on the medians of $ G_{0} $ and $ G_{5} $:
\[ 
H_{0}: (\xi_{0}, \xi_{5}) = (\xi_{0}^{*}, \xi_{5}^{*}) \, \, \, \text { versus } 
H_{1}: (\xi_{0}, \xi_{5}) \neq (\xi_{0}^{*}, \xi_{5}^{*})
\]
where $ \xi_{0}^{*}, \xi_{5}^{*} $ are the true values. 
Note that although we simulate data from normal distributions, the parametric information 
does not play any role in the data analysis.

\begin{figure}[!ht] 
\centering
\caption{Q-Q plot of $ R_{n} $ values against $ \chi_{2}^{2} $ based on 
normal data of equal sample size $ n_{r} = 100 $.}
\label{Normal_QQ}
\includegraphics[height = 241.5pt, width = 241.5pt]{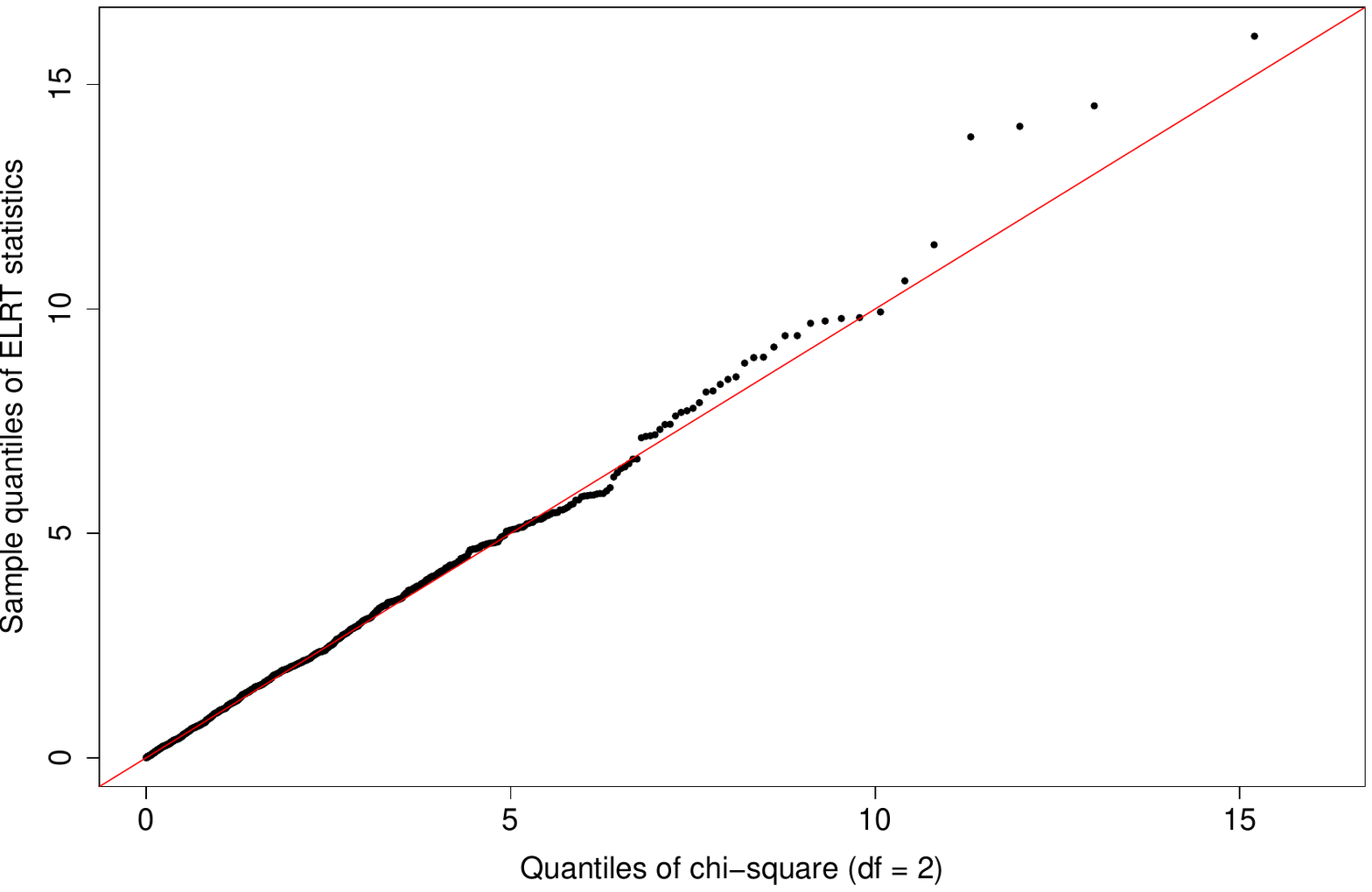}
\end{figure}

Because $ H_{0} $ is true, $ R_{n} $ has a $ \chi_{2}^{2} $ limiting distribution.
Figure~\ref {Normal_QQ} gives a quantile-quantile (Q-Q) plot of the 1000 simulated 
$ R_{n} $ values against the $ \chi_{2}^{2} $ distribution. 
Over the range from 0 to 6 that matters in most applications, the points are close to the red 45-degree line.
Clearly, the chi-square distribution is a good approximation of the sampling distribution of $ R_{n} $, 
demonstrating good agreement with Theorem \ref {thm:ELRT_Stat_chisq}.   

\begin{figure}[!ht] 
\begin{center}
\caption{Confidence regions of $ (\xi_{0}, \xi_{5}) $ by ELRT (solid), Wald (dashed), and nonparametric (dotted) methods, 
based on a simulated normal data set of equal sample size $ n_{r} = 100 $. 
The true quantiles are marked with a diamond. The level of confidence is $ 95\% $.}
\label{CIs_Normal1}
\includegraphics[height = 8cm, width = 345.0pt]{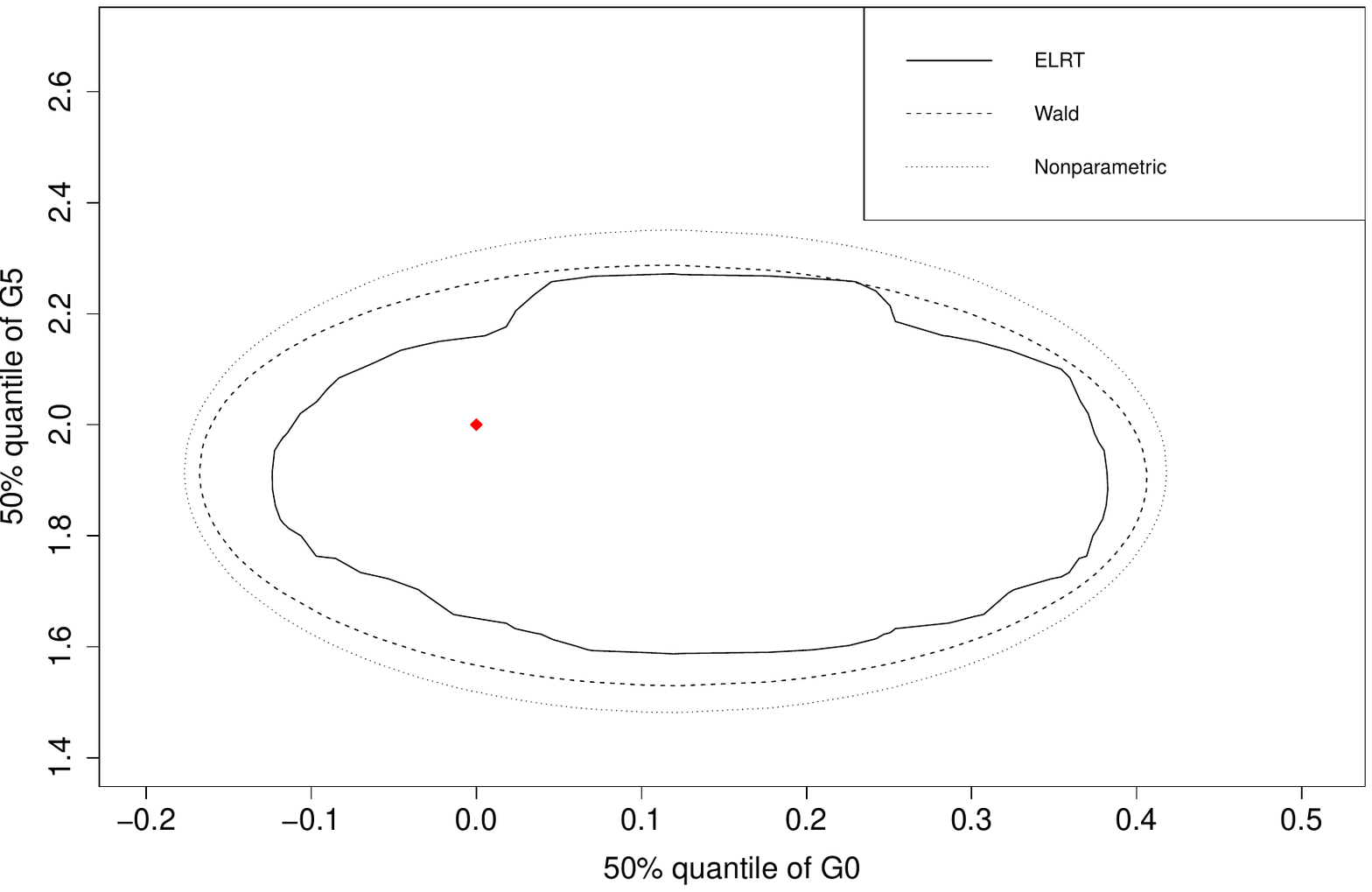}
\end{center}
\end{figure}

In Figure~\ref {CIs_Normal1}, we depict the $ 95\% $ confidence regions of 
$ \bxi = (\xi_{0}, \xi_{5}) $ based on the ELRT in \eqref {ELRT_CR}, 
the Wald method in \eqref {Wald_CR}, and the nonparametric method in \eqref {Nonpara_CR} 
based on a typical simulated data set with the true $ \bxi^{*} $ marked as a red diamond. 
The ELRT contour is not smooth because $ R_{n} (\bxi) $ is not smooth at data points.  
Clearly, the ELRT confidence region has the smallest area and is therefore the most efficient. 
In Table~\ref {Normal_Ava_Ecv1}, we make direct quantitative comparisons 
between the three methods in terms of the coverage probabilities and 
areas of the $ 90\% $ and $ 95\% $ confidence regions.
We remark that the ELRT confidence region can be approximated
by triangles all pointing to the MELE.
We add up the areas of these triangles to get the total area.
Both the LRT and Wald methods under the DRM have empirical coverage probabilities 
close to the nominal levels; the nonparametric method has overcoverage.
In general, the ELRT outperforms.

\begin{table}[!ht]
\centering
\caption{Empirical coverage probabilities and average areas based on normal data of equal sample size.}
\label{Normal_Ava_Ecv1}
\resizebox{345.0pt}{!}{
\begin{tabular}{cccccc}
\toprule
\multirow{2}{*}{} & \multirow{2}{*}{Method} & \multicolumn{2}{c}{90\%} & \multicolumn{2}{c}{95\%}\\
\cmidrule(lr){3-4} \cmidrule(lr){5-6} 
& & Coverage probability & Area & Coverage probability & Area \\
\midrule
\multicolumn{6}{c}{\multirow{2}{*}{$ n_{r} = 100 $}} \\ 
\\ 
& ELRT & $ 89.1\% $ & $ 0.250 $ & $ 95.8\% $ & $ 0.323 $  \\ 
& Wald & $ 90.8\% $ & $ 0.266 $ & $ 95.4\% $ & $ 0.347 $  \\
& Nonparametric & $ 91.7\% $ & $ 0.374 $ & $ 95.9\% $ & $ 0.487 $  \\
\multicolumn{6}{c}{\multirow{2}{*}{$ n_{r} = 200 $}} \\ 
\\ 
& ELRT & $ 89.7\% $ & $ 0.126 $ & $ 95.0\% $ & $ 0.164 $  \\ 
& Wald & $ 90.5\% $ & $ 0.132 $ & $ 95.2\% $ & $ 0.171 $  \\
& Nonparametric & $ 90.3\% $ & $ 0.183 $ & $ 95.3\% $ & $ 0.239 $  \\
\bottomrule
\end{tabular}
}
\end{table}

\begin{figure}[!ht] 
\begin{center}
\caption{Confidence regions of $ (\xi_{2}, \xi_{3}) $ by ELRT (solid), Wald (dashed), and nonparametric (dotted) methods, 
based on a simulated normal data set of unequal sample sizes. 
The true quantiles are marked with a diamond. The level of confidence is $ 95\% $.}
\label{CIs_Normal2}
\includegraphics[height = 8cm, width = 345.0pt]{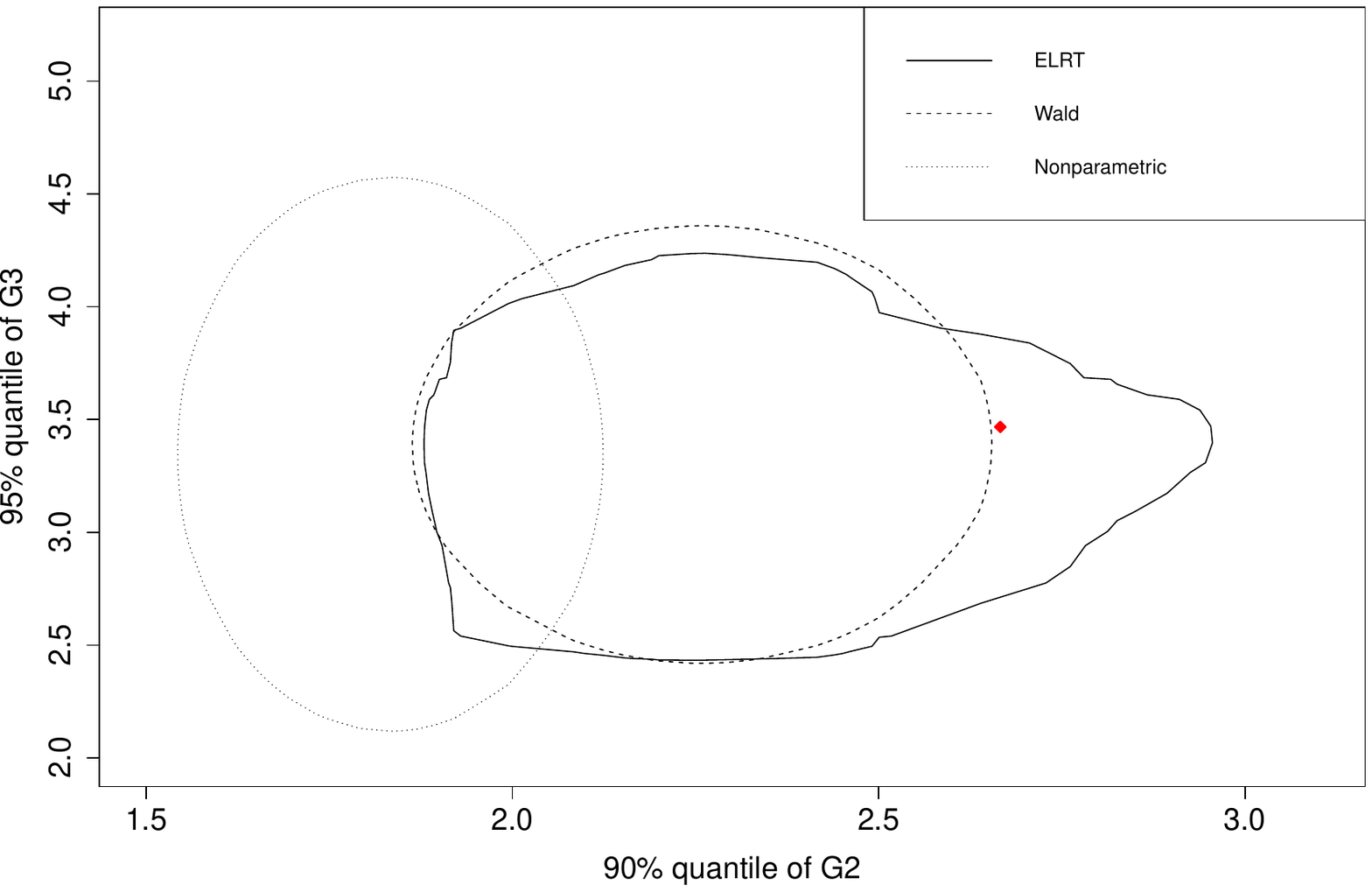}
\end{center}
\end{figure}

In applications, the sample sizes from different populations are unlikely to be equal.
Does the superiority of the ELRT require equal sample sizes from these populations? 
We also simulated data from the same distributions with unequal sample sizes. 
We set the sizes of populations $ G_{0}, G_{1}, G_{4}, G_{5} $ to 100 and 200, 
and the sizes of populations $ G_{2}, G_{3} $ to 50 and 100, respectively.  
We constructed confidence regions for the $ 90$th percentile of $ G_{2} $ and the $ 95 $th 
percentile of $ G_{3} $, where both populations have the smaller sample sizes.  
Figure~\ref {CIs_Normal2} shows the three $ 95\% $ confidence regions based on a simulated data set; we see that the ELRT is superior.
Admittedly, this is one of the more extreme cases. 
Table~\ref {Normal_Ava_Ecv2} gives the average areas and empirical coverage probabilities of the three confidence regions, based on 1000 repetitions. 
The ELRT confidence regions have the most accurate coverage probabilities, while the other two methods have low coverage. 
The ELRT confidence regions have larger average areas that are not excessive.

\begin{table}[!ht]
\centering
\caption{Empirical coverage probabilities and average areas based on normal data of unequal sample sizes.}
\label{Normal_Ava_Ecv2}
\resizebox{345.0pt}{!}{
\begin{tabular}{cccccc}
\toprule
\multirow{2}{*}{} & \multirow{2}{*}{Method} & \multicolumn{2}{c}{90\%} & \multicolumn{2}{c}{95\%}\\
\cmidrule(lr){3-4} \cmidrule(lr){5-6} 
& & Coverage probability & Area & Coverage probability & Area \\
\midrule
\multicolumn{6}{c}{\multirow{2}{*}{$ n_{2} =  n_{3} = 50, n_{0} = n_{1} = n_{4} = n_{5} = 100 $}} \\ 
\\ 
& ELRT & $ 90.1\% $ & $ 1.307 $ & $ 94.5\% $ & $ 1.741 $  \\ 
& Wald & $ 83.7\% $ & $ 1.096 $ & $ 88.9\% $ & $ 1.427 $  \\
& Nonparametric & $ 73.6\% $ & $ 1.439 $ & $ 80.0\% $ & $ 1.873 $  \\
\multicolumn{6}{c}{\multirow{2}{*}{$ n_{2} =  n_{3} = 100, n_{0} = n_{1} = n_{4} = n_{5} = 200 $}} \\ 
\\ 
& ELRT & $ 90.1\% $ & $ 0.642 $ & $ 94.5\% $ & $ 0.843 $  \\ 
& Wald & $ 86.7\% $ & $ 0.572 $ & $ 91.8\% $ & $ 0.744 $  \\
& Nonparametric & $ 81.3\% $ & $ 0.804 $ & $ 86.7\% $ & $ 1.046 $  \\
\bottomrule
\end{tabular}
}
\end{table}

\subsection{Data generated from gamma distributions} 

In applications, income, lifetime, expenditure, and strength data are positive and skewed.
Gamma or Weibull distributions are often used for statistical inference in such applications. 
In the presence of multiple samples, replacing the parametric model by a DRM 
with $ \bq (x) = (1, x, \log x)^{\top} $ is an attractive option to reduce the risk of model mis-specification. 
We generate 1000 sets of $ m+1 = 6 $ independent samples of sizes $ n_{r} = 100 \text { and } 200 $ from gamma distributions with shape parameters $ (5, 5, 6, 6, 7, 7) $ and scale parameters $ (2, 1.9, 1.8, 1.7, 1.6, 1.5) $.
We test the hypothesis on the medians of $ G_{1} $ and $ G_{2} $: 
\[ 
H_{0}: (\xi_{1}, \xi_{2}) = (\xi_{1}^{*}, \xi_{2}^{*}) \, \, \, 
\text { versus } H_{1}: (\xi_{1}, \xi_{2}) \neq (\xi_{1}^{*}, \xi_{2}^{*}), 
\] 
where $ \xi_{1}^{*}, \xi_{2}^{*} $ are the true medians of 
$ \mathrm {Gamma} (5, 1.9) $ and $ \mathrm {Gamma} (6, 1.8) $, respectively. 
Note that although we simulate data from gamma distributions, 
the parametric information does not play any role in the data analysis.

\begin{figure}[!ht] 
\centering
\caption{Q-Q plot of $ R_{n} $ values against $ \chi_{2}^{2} $ based on gamma data 
of equal sample size $ n_{r} = 100 $.}
\label{Gamma_QQ}
\includegraphics[height = 241.5pt, width = 241.5pt]{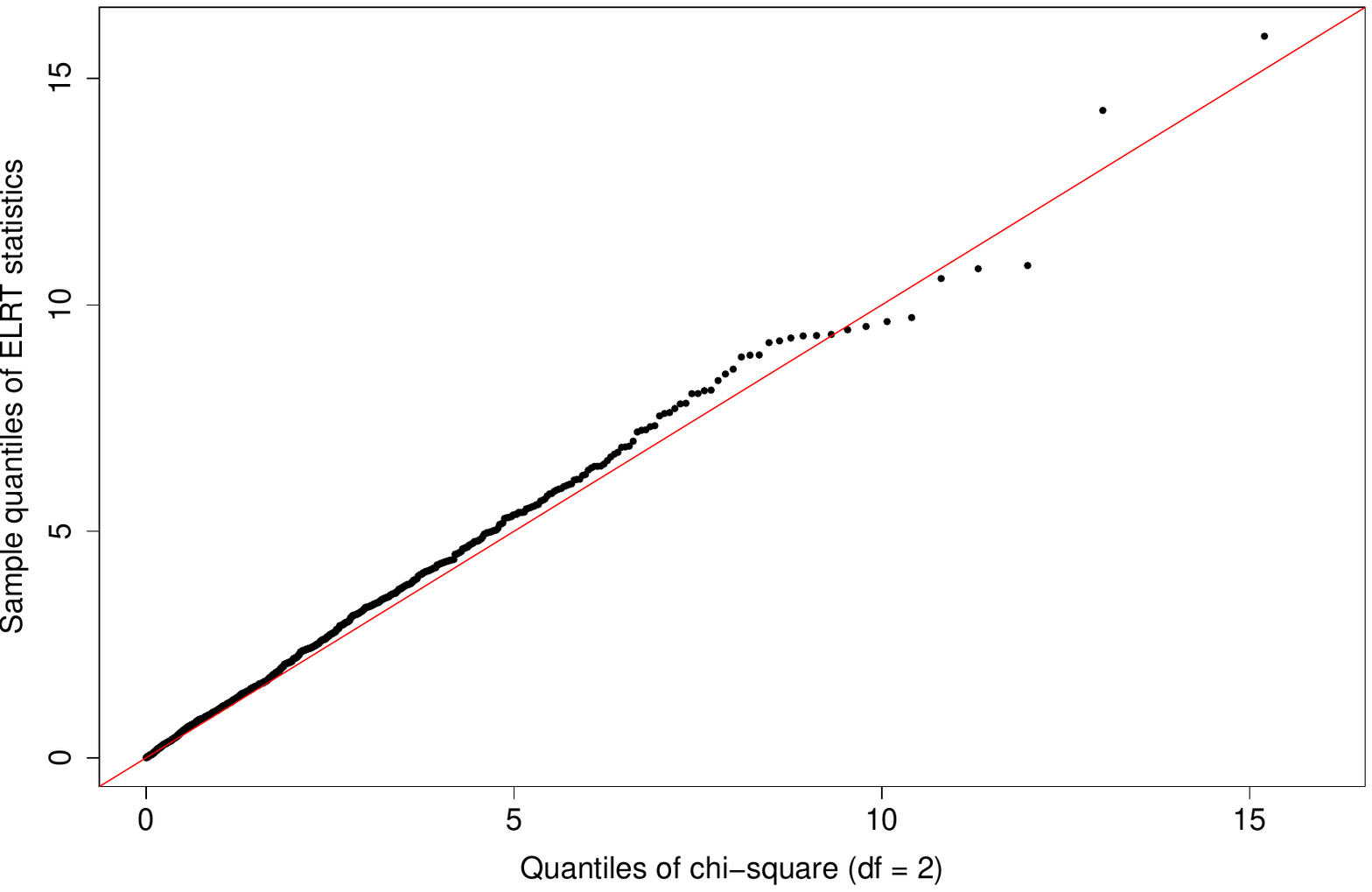}
\end{figure}

Figure~\ref {Gamma_QQ} shows the Q-Q plot based on 1000 $ R_{n} $ 
values against the theoretical limiting distribution $ \chi_{2}^{2} $. 
The points in the Q-Q plot are close to (but slightly above) the 45-degree line in the range from 0 to 6. 
This implies that the corresponding tests will have close to nominal levels.
Overall, the chi-square approximation is satisfactory.

\begin{figure}[!ht] 
\centering
\caption{Confidence regions of $ (\xi_{1}, \xi_{2}) $ by ELRT (solid), Wald (dashed), 
and nonparametric (dotted) methods, 
based on a simulated gamma data set of equal sample size $ n_{r} = 100 $. 
The true quantiles are marked with a diamond. The level of confidence is $ 95\% $.}
\label{Gamma_CIs1}
\includegraphics[height = 8cm, width = 345.0pt]{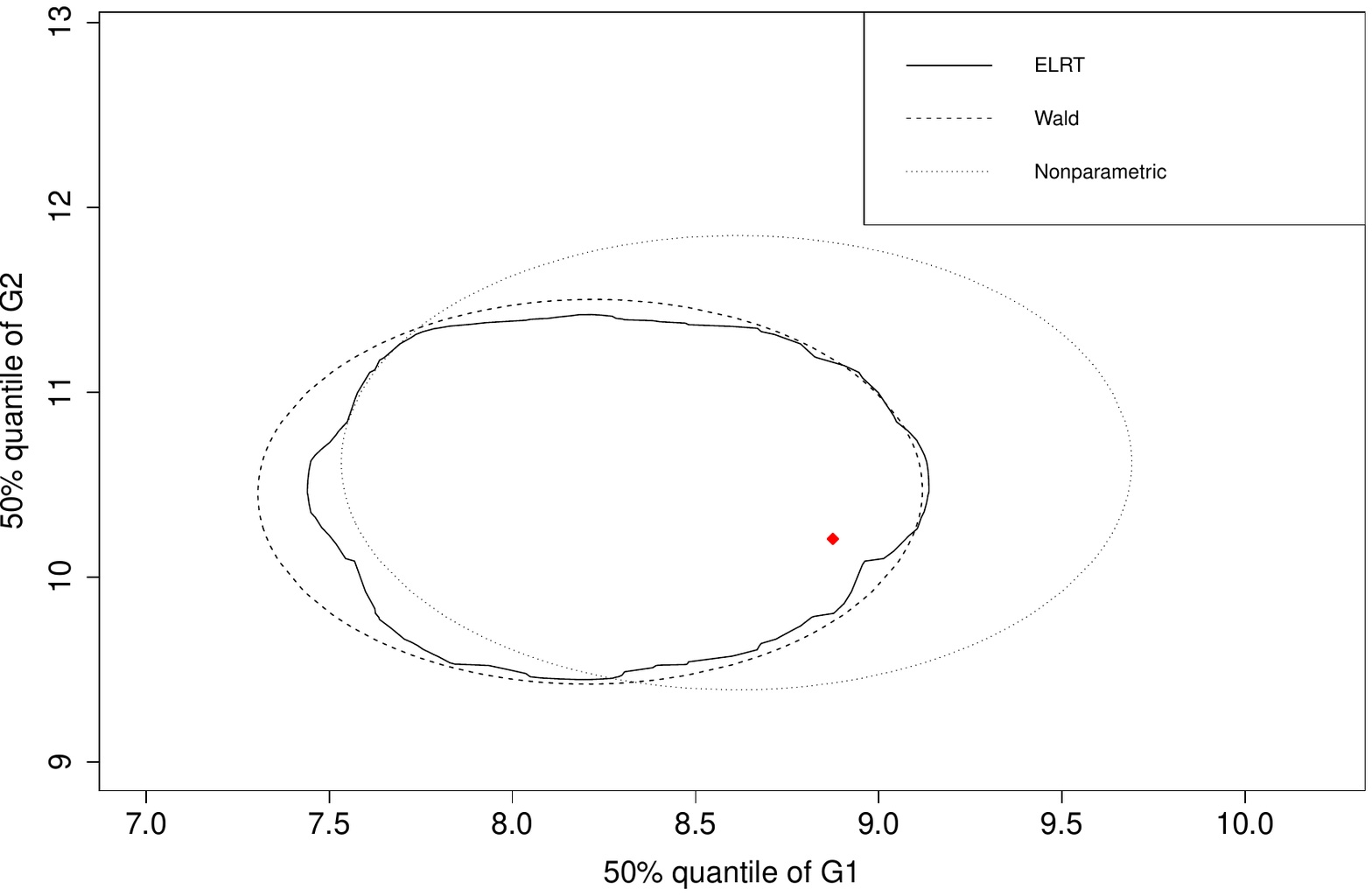}
\end{figure}

In Figure~\ref {Gamma_CIs1}, we depict the $ 95\% $ confidence regions of 
$ \bxi = (\xi_{1}, \xi_{2}) $ using the ELRT in \eqref {ELRT_CR}, 
the Wald method in \eqref {Wald_CR}, and the nonparametric method in 
\eqref {Nonpara_CR}, based on a typical simulated data set with $ \bxi^{*} $ marked as a red diamond. 
Clearly, the ELRT-based confidence region has a smaller area and is therefore more efficient. 
In Table~\ref {Gamma_Ava_Ecv1} we make direct quantitative comparisons 
of the coverage probabilities and areas.
Both the ELRT and Wald methods under the DRM have empirical 
coverage probabilities very close to the nominal levels.
The nonparametric confidence regions have overcoverage and inflated sizes.
We again conclude that the ELRT is superior to the nonparametric method. 

\begin{table}[!ht]
\centering
\caption{Empirical coverage probabilities and average areas based on gamma data of equal sample size.}
\label{Gamma_Ava_Ecv1}
\resizebox{345.0pt}{!}{
\begin{tabular}{cccccc}
\toprule
\multirow{2}{*}{} & \multirow{2}{*}{Method} & \multicolumn{2}{c}{90\%} & \multicolumn{2}{c}{95\%}\\
\cmidrule(lr){3-4} \cmidrule(lr){5-6} 
& & Coverage probability & Area & Coverage probability & Area \\
\midrule
\multicolumn{6}{c}{\multirow{2}{*}{$ n_{r} = 100 $}} \\ 
\\ 
& ELRT & $ 88.3\% $ & $ 2.808 $ & $ 94.2\% $ & $ 3.665 $  \\ 
& Wald & $ 89.9\% $ & $ 2.953 $ & $ 95.3\% $ & $ 3.843 $  \\
& Nonparametric & $ 92.1\% $ & $ 4.264 $ & $ 95.2\% $ & $ 5.547 $  \\
\multicolumn{6}{c}{\multirow{2}{*}{$ n_{r} = 200 $}} \\ 
\\ 
& ELRT & $ 88.6\% $ & $ 1.395 $ & $ 94.4\% $ & $ 1.822 $  \\ 
& Wald & $ 89.7\% $ & $ 1.451 $ & $ 95.3\% $ & $ 1.889 $  \\
& Nonparametric & $ 89.3\% $ & $ 2.111 $ & $ 94.3\% $ & $ 2.747 $  \\
\bottomrule
\end{tabular}
}
\end{table}

\begin{figure}[!ht] 
\centering
\caption{Confidence regions of $ (\xi_{4}, \xi_{5}) $ by ELRT (solid), Wald (dashed), and nonparametric (dotted) methods, 
based on a simulated gamma data set of equal sample size $ n_{r} = 100 $. 
The true quantiles are marked with a diamond. The level of confidence is $ 95\% $.}
\label{Gamma_CIs2}
\includegraphics[height = 8cm, width = 345.0pt]{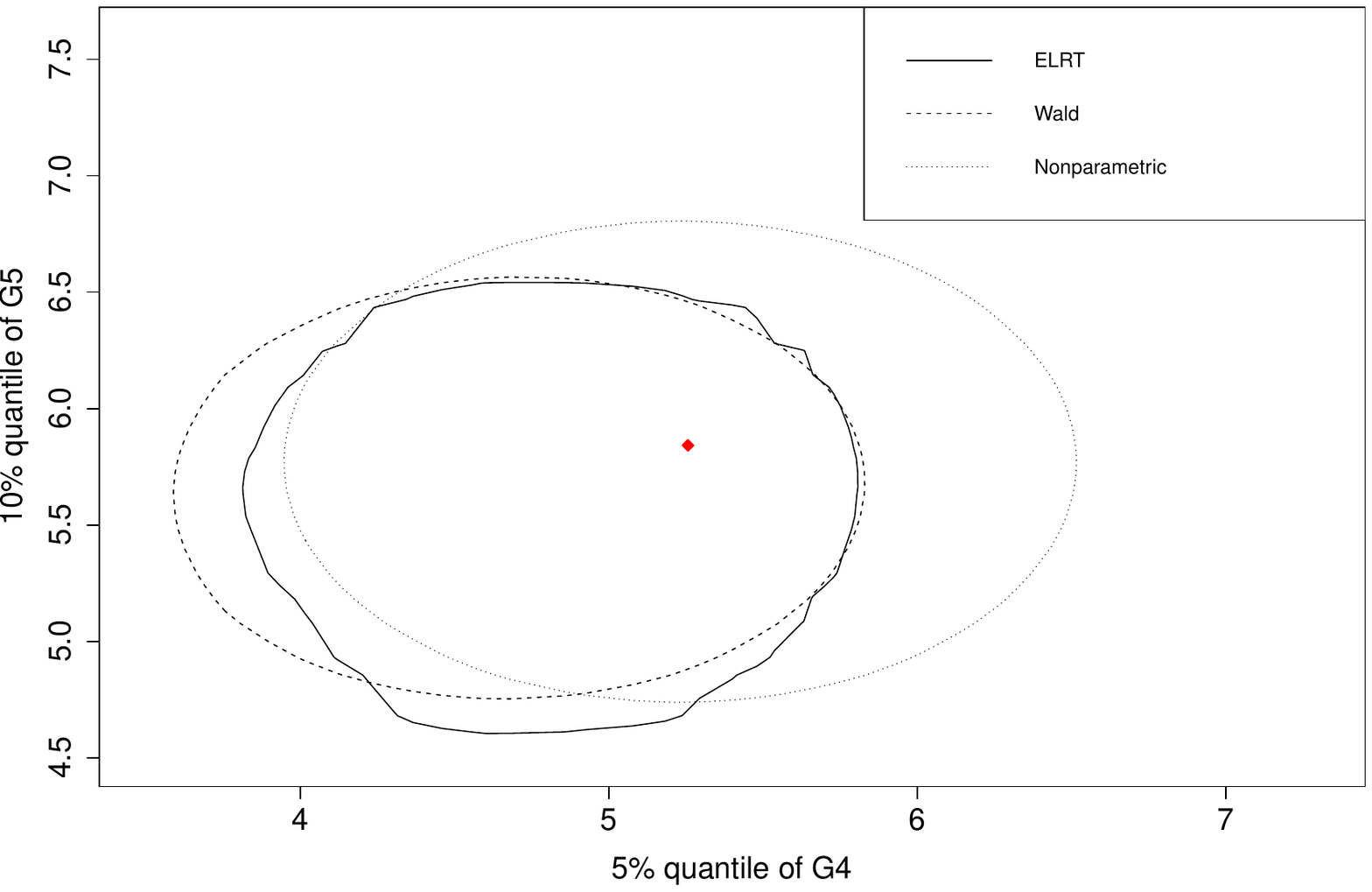}
\end{figure}

We also study the confidence regions for a pair of lower quantiles: the $ 5 $th 
percentile of $ G_{4} $ and the $ 10 $th percentile of $ G_{5} $. 
Figure~\ref {Gamma_CIs2} shows the three $ 95\% $ confidence regions based on a simulated data set. 
Table~\ref {Gamma_Ava_Ecv2} gives the average areas and coverage probabilities of the three confidence regions, based on 1000 repetitions. 
The ELRT method is still the most efficient. 
Maintaining the accurate coverage probabilities, the ELRT confidence regions still have satisfactory areas that are comparable to the Wald confidence regions.

\begin{table}[!ht]
\centering
\caption{Empirical coverage probabilities and average areas based on gamma data of equal sample size.}
\label{Gamma_Ava_Ecv2}
\resizebox{345.0pt}{!}{
\begin{tabular}{cccccc}
\toprule
\multirow{2}{*}{} & \multirow{2}{*}{Method} & \multicolumn{2}{c}{90\%} & \multicolumn{2}{c}{95\%}\\
\cmidrule(lr){3-4} \cmidrule(lr){5-6} 
& & Coverage probability & Area & Coverage probability & Area \\
\midrule
\multicolumn{6}{c}{\multirow{2}{*}{$ n_{r} = 100 $}} \\ 
\\ 
& ELRT & $ 88.4\% $ & $ 2.312 $ & $ 93.7\% $ & $ 3.022 $  \\ 
& Wald & $ 86.5\% $ & $ 2.236 $ & $ 92.0\% $ & $ 2.910 $  \\
& Nonparametric & $ 82.8\% $ & $ 3.250 $ & $ 88.7\% $ & $ 4.229 $  \\
\multicolumn{6}{c}{\multirow{2}{*}{$ n_{r} = 200 $}} \\ 
\\ 
& ELRT & $ 90.8\% $ & $ 1.139 $ & $ 95.3\% $ & $ 1.486 $  \\ 
& Wald & $ 90.4\% $ & $ 1.114 $ & $ 95.4\% $ & $ 1.449 $  \\
& Nonparametric & $ 87.0\% $ & $ 1.684 $ & $ 92.6\% $ & $ 2.191 $  \\
\bottomrule
\end{tabular}
}
\end{table}

%%%%%%%%%%%%%%%%%%%%%%%%%%%%%%%%%%%%%%%%%%%%%%%%%%%%%%%%%%%%%%%%%%%%%%

\section{Real-data analysis} 
\label{sec:realdata}

In the previous simulations, we chose the most suitable basis function $ \bq (x) $ in each case because the population distributions were known to us.
This is not possible in real-world applications.
In this section, we create a simulation population based on the US Consumer Expenditure Surveys data concerning US expenditure, income, and demographics.
The data set is available on the US Bureau of Labor Statistics website (\url {https://www.bls.gov/cex/pumd.htm}). 
The data are collected by the Census Bureau in the form of panel surveys, in which approximately 5000 households are contacted each quarter. 
After a household has been surveyed it is dropped from subsequent surveys and replaced by a new household. 
The response variable is the annual sum of the wages or salary income received by all household members before any deductions.
Household income is a good reflection of economic well-being. 
The data files include some imputed values to replace missing values due to non-response.

\begin{figure}[!ht] 
\centering
\caption{Histograms of log-transformed annual household incomes.} 
\label{RealData_hist}
\includegraphics[height = 10cm, width = 345.0pt]{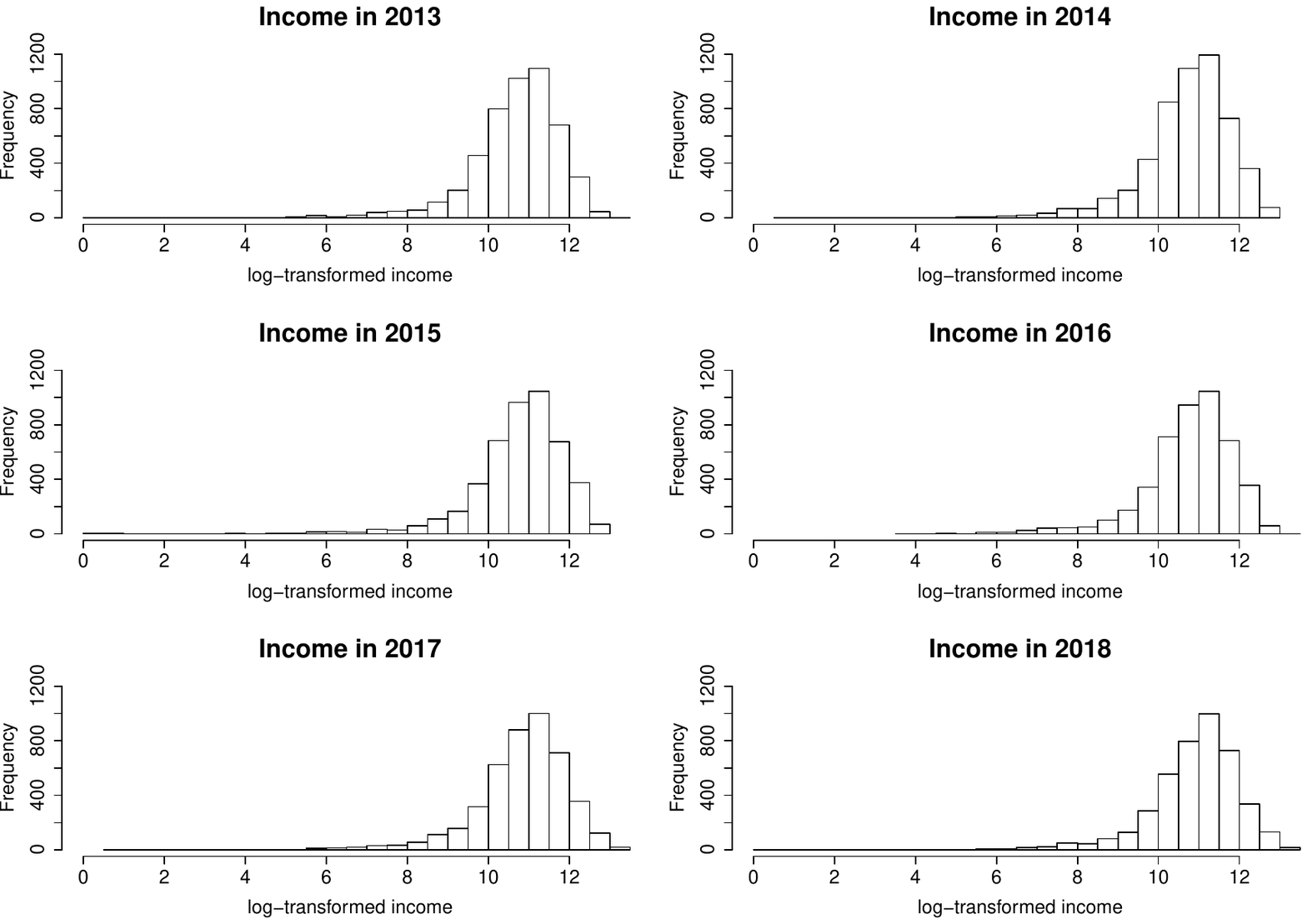}
\end{figure}

We study a six-year period from 2013 to 2018, 
and we log-transform the response values to make the scale more suitable for numerical computation. 
Note that the quantiles are transformation equivariant.
We exclude households that have no recorded income even after imputation, 
and there remain 4919, 5304, 4641, 4606, 4475, and 4222 households from 2013 to 2018.
The histograms shown in Figure~\ref {RealData_hist} indicate that 
it is difficult to determine a suitable parametric model for these data sets, but a DRM may work well enough. 
We take the basis function $ \bq (x) = (1, x, x^{2})^{\top} $; it may not be the best choice, but as a result the
simulation results for the DRM analysis are more convincing. 

\begin{figure}[!ht] 
\centering
\caption{Q-Q plots of $ R_{n} $ values against $ \chi_{1}^{2} $, 
based on real data of equal sample size $ n_{r} = 100 $. Quantile levels are $ 20\% $.}
\label{RealData_QQ_ELRT_20quan}
\includegraphics[height = 6cm, width = 345.0pt]{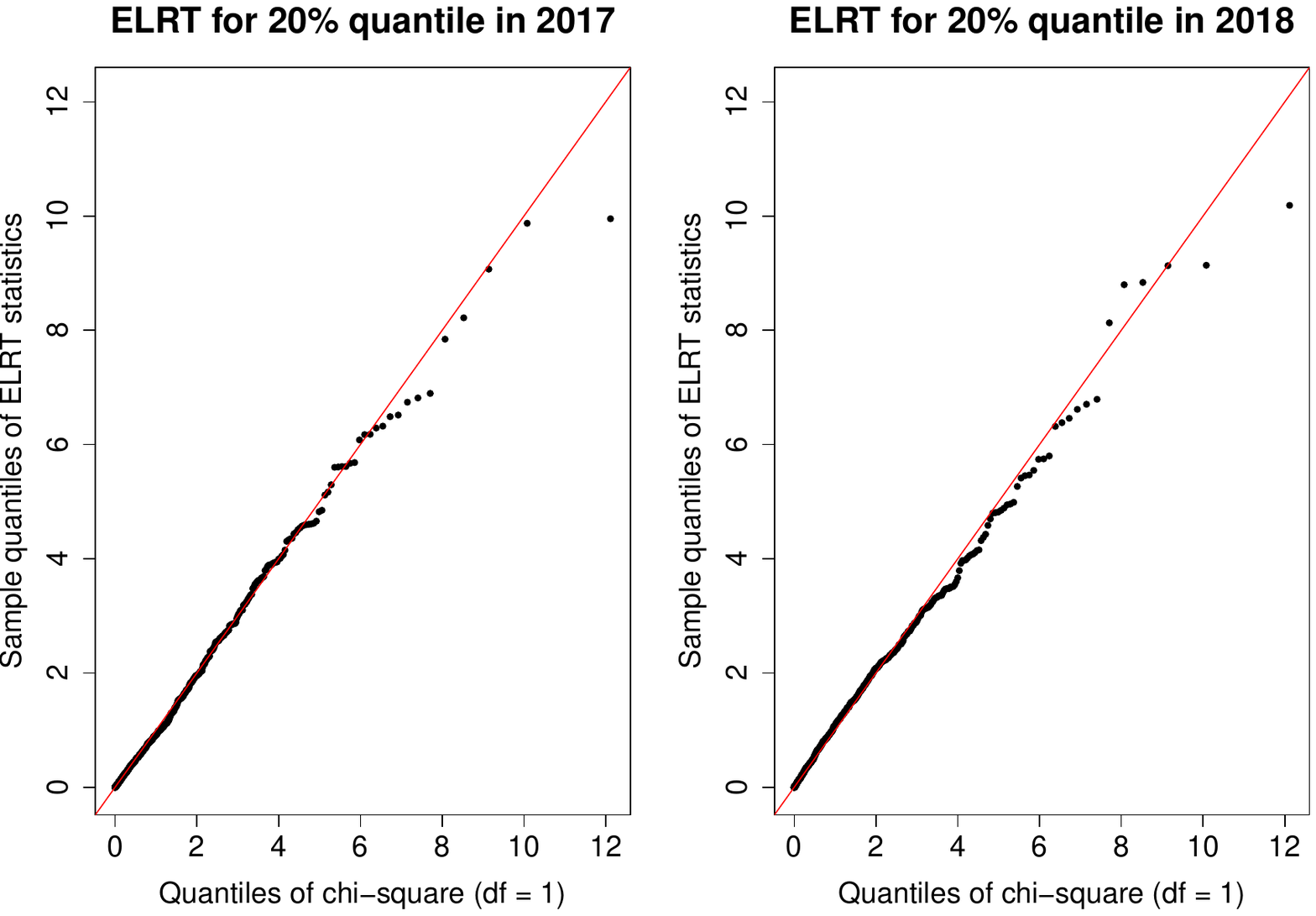} 
\end{figure}

\begin{figure}[!ht] 
\centering
\caption{Q-Q plots of $ R_{n} $ values against $ \chi_{1}^{2} $, 
based on real data of equal sample size $ n_{r} = 100 $. Quantile levels are $ 50\% $.}
\label{RealData_QQ_ELRT_50quan}
\includegraphics[height = 6cm, width = 345.0pt]{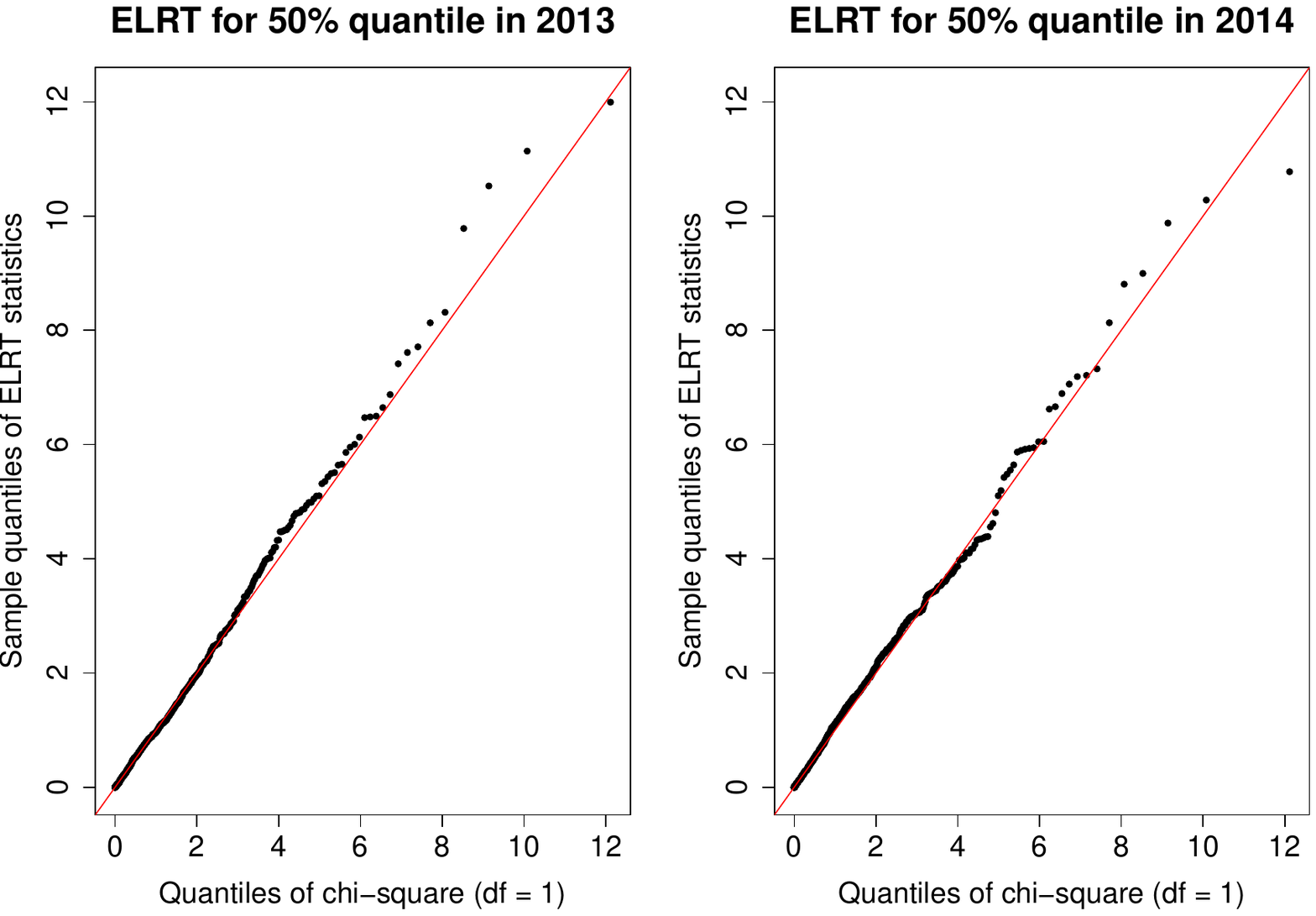} 
\end{figure}

In this simulation, we form 6 populations based on the yearly data sets.
We test hypotheses on the $ 20$th and $ 50$th percentiles 
based on independent samples of size 100, which are obtained by 
sampling with replacement from the respective populations.
To test the size of a single quantile of a single population, the limiting distribution of 
$ R_{n} $ is $ \chi_{1}^{2} $.
Figures~\ref {RealData_QQ_ELRT_20quan} and \ref {RealData_QQ_ELRT_50quan} 
contain a few Q-Q plots of $ R_{n} $ versus $ \chi_{1}^{2} $ 
for $ H_{0}: \xi_{r} = \xi_{r}^{*} $ with $ \tau_{r} = 20\% $ or $ \tau_{r} = 50\% $. 
In all the plots, the points of $ R_{n} $ are close to the 45-degree line.
Thus, the precision of the chi-square approximation is satisfactory.
The plots for other levels or populations are similar and not presented.

\begin{figure}[!ht] 
\centering
\caption{Q-Q plots of Wald statistic values against $ \chi_{1}^{2} $, 
based on real data of equal sample size $ n_{r} = 100 $. Quantile levels are $ 20\% $.}
\label{RealData_QQ_Wald_20quan}
\includegraphics[height = 6cm, width = 345.0pt]{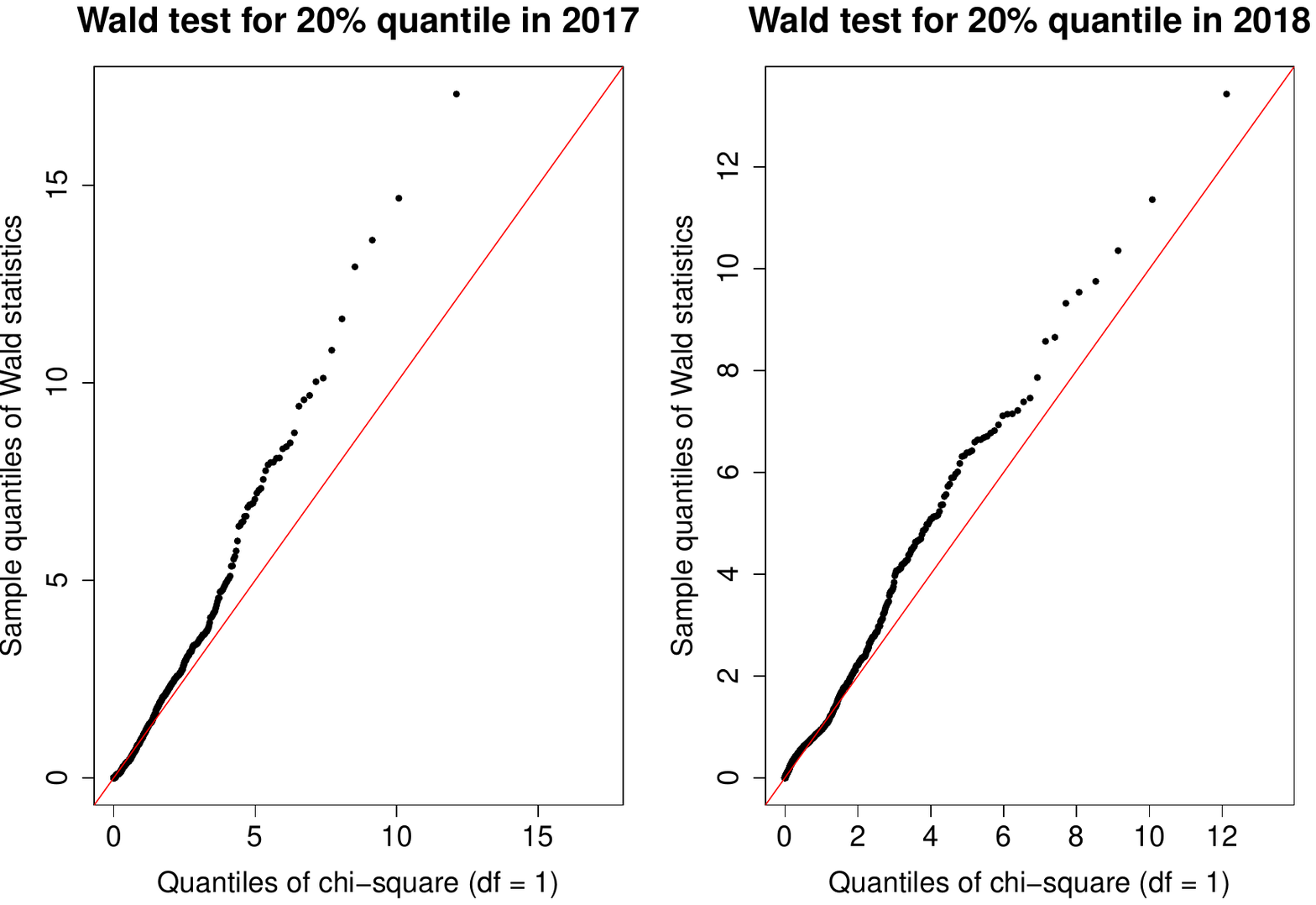} 
\end{figure}

\begin{figure}[!ht] 
\centering
\caption{Q-Q plots of Wald statistic values against $ \chi_{1}^{2} $, 
based on real data of equal sample size $ n_{r} = 100 $. Quantile levels are $ 50\% $.}
\label{RealData_QQ_Wald_50quan}
\includegraphics[height = 6cm, width = 345.0pt]{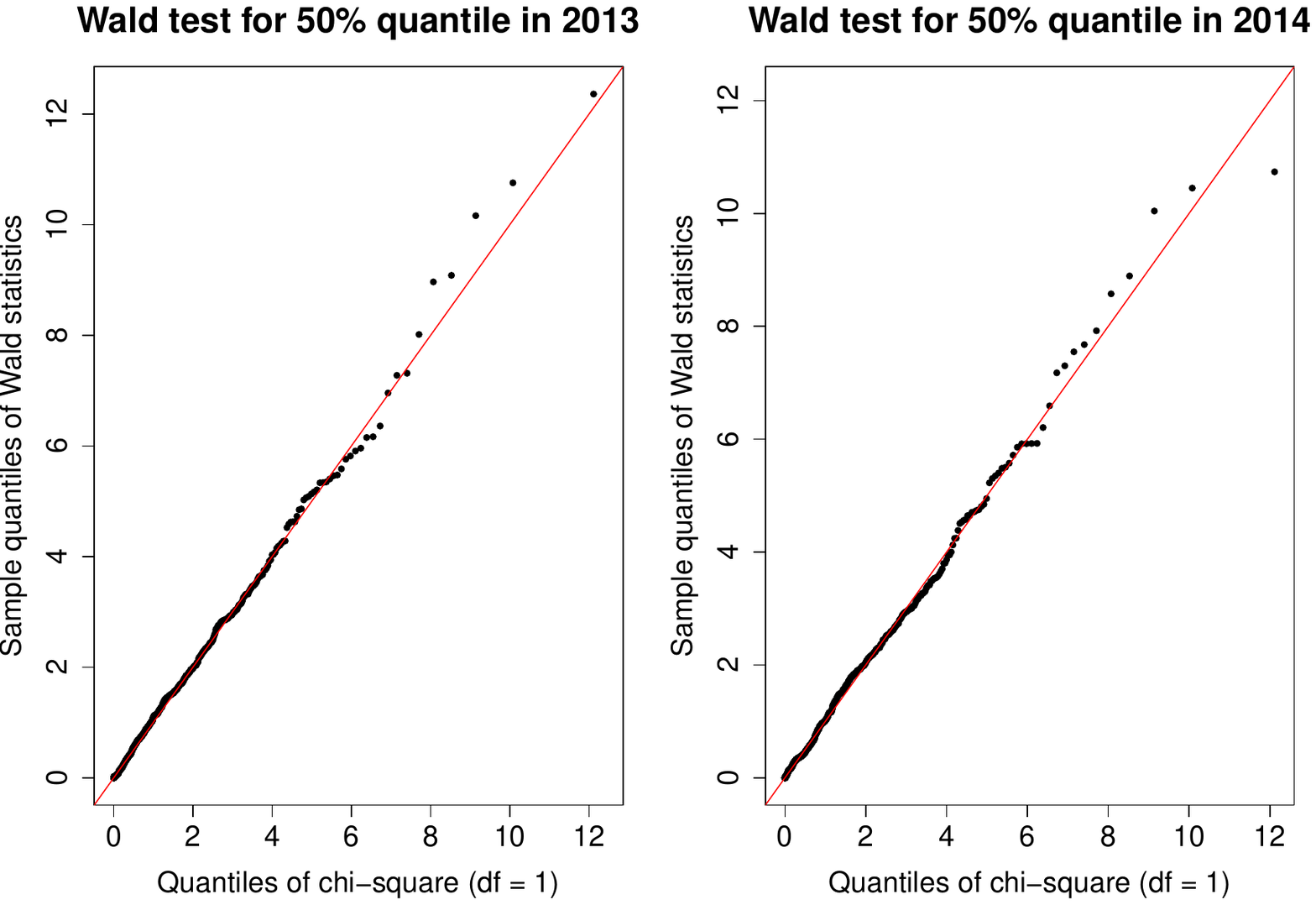} 
\end{figure}

The Wald method \eqref {Wald_CR} may be regarded as being derived from an asymptotic $ \chi_{1}^{2} $ distributed statistic: 
\[
W_{n} = n (\tilde \xi_{r} - \xi_{r}^{*})^{\top} \tilde \Omega^{-1} (\tilde \xi_{r} - \xi_{r}^{*}).
\]
We also obtain $ W_{n} $ values and construct Q-Q plots, and a selected few are given in Figures~\ref {RealData_QQ_Wald_20quan} and \ref {RealData_QQ_Wald_50quan}.
These plots show that the chi-square approximation is not as satisfactory.
There are many possible explanations, but a major factor could be the unstable variance estimator $ \tilde \Omega $ that the Wald method must rely on, especially for lower quantiles. 
One of the most valued properties of the likelihood ratio test approach is that there is no need to estimate a scale factor. 

A direct consequence of the poor chi-square approximation could be undercoverage of the confidence intervals. 
Table~\ref {RealData_CI_Ava_Ecv1} gives the coverage probabilities and average lengths of the confidence intervals based on three methods: ELRT in \eqref {ELRT_CR}, Wald in \eqref {Wald_CR}, and nonparametric in \eqref {Nonpara_CR}. 
The improved efficiency of the DRM is best reflected in the average lengths of the confidence intervals. 
It can be seen that the DRM-based methods achieve on average about $ 15\% $ and $ 25\% $ improvement over the nonparametric method for the $ 20 $th and $ 50 $th percentiles respectively. 
Comparing the ELRT and Wald methods, both done under DRM, we find that the ELRT is comparable to the Wald method for the $ 20 $th percentile and clearly more efficient for the $ 50 $th percentile.

\begin{table}[!ht]
\centering
\caption{Average lengths and empirical coverage probabilities of the individual confidence intervals, based on real data of equal sample size $ n_{r} = 100 $.}
\label{RealData_CI_Ava_Ecv1}
%\resizebox{345.0pt}{!}{
\begin{tabular}{cccccccc}
\toprule
\multirow{2}{*}{} & \multirow{2}{*}{Year} & \multicolumn{2}{c}{ELRT} & \multicolumn{2}{c}{Wald} & \multicolumn{2}{c}{Nonparametric} \\
\cmidrule(lr){3-4} \cmidrule(lr){5-6} \cmidrule(lr){7-8} 
& & $ 90\% $ & $ 95\% $ & $ 90\% $ & $ 95\% $ & $ 90\% $ & $ 95\% $ \\
\midrule
\multicolumn{8}{c}{\multirow{1}{*}{{\bf Average lengths}}} \\ 
\midrule 
\multicolumn{8}{c}{\multirow{2}{*}{quantile levels all $ = 20\% $}} \\ 
\\ 
& 2013 & $ 0.465 $ & $ 0.563 $ & $ 0.440 $ & $ 0.524 $ & $ 0.513 $ & $ 0.611 $ \\
& 2014 & $ 0.464 $ & $ 0.559 $ & $ 0.437 $ & $ 0.520 $ & $ 0.528 $ & $ 0.630 $ \\
& 2015 & $ 0.459 $ & $ 0.553 $ & $ 0.432 $ & $ 0.515 $ & $ 0.519 $ & $ 0.619 $ \\ 
& 2016 & $ 0.461 $ & $ 0.558 $ & $ 0.435 $ & $ 0.519 $ & $ 0.527 $ & $ 0.628 $ \\ 
& 2017 & $ 0.459 $ & $ 0.557 $ & $ 0.434 $ & $ 0.518 $ & $ 0.539 $ & $ 0.642 $ \\ 
& 2018 & $ 0.438 $ & $ 0.529 $ & $ 0.416 $ & $ 0.496 $ & $ 0.523 $ & $ 0.623 $ \\
& average & $ 0.458 $ & $ 0.553 $ & $ 0.433 $ & $ 0.515 $ & $ 0.525 $ & $ 0.626 $ \\ 
\multicolumn{8}{c}{\multirow{2}{*}{quantile levels all $ = 50\% $}} \\ 
\\ 
& 2013 & $ 0.307 $ & $ 0.364 $ & $ 0.315 $ & $ 0.376 $ & $ 0.383 $ & $ 0.457 $ \\
& 2014 & $ 0.306 $ & $ 0.366 $ & $ 0.316 $ & $ 0.376 $ & $ 0.379 $ & $ 0.452 $ \\
& 2015 & $ 0.304 $ & $ 0.364 $ & $ 0.314 $ & $ 0.374 $ & $ 0.374 $ & $ 0.446 $ \\ 
& 2016 & $ 0.305 $ & $ 0.364 $ & $ 0.315 $ & $ 0.375 $ & $ 0.382 $ & $ 0.455 $ \\ 
& 2017 & $ 0.304 $ & $ 0.364 $ & $ 0.316 $ & $ 0.376 $ & $ 0.390 $ & $ 0.465 $ \\ 
& 2018 & $ 0.300 $ & $ 0.357 $ & $ 0.311 $ & $ 0.371 $ & $ 0.373 $ & $ 0.444 $ \\
& average & $ 0.304 $ & $ 0.363 $ & $ 0.315 $ & $ 0.375 $ & $ 0.380 $ & $ 0.453 $ \\ 
\\ 
\midrule 
\multicolumn{8}{c}{\multirow{1}{*}{{\bf Empirical coverage probabilities}}} \\ 
\midrule
\multicolumn{8}{c}{\multirow{2}{*}{quantile levels all $ = 20\% $}} \\ 
\\ 
& 2013 & $ 88.0\% $ & $ 94.0\% $ & $ 88.7\% $ & $ 93.2\% $ & $ 87.7\% $ & $ 92.3\% $ \\
& 2014 & $ 90.1\% $ & $ 95.1\% $ & $ 88.7\% $ & $ 94.7\% $ & $ 87.9\% $ & $ 92.6\% $ \\
& 2015 & $ 89.8\% $ & $ 94.6\% $ & $ 88.6\% $ & $ 93.6\% $ & $ 89.5\% $ & $ 94.3\% $ \\ 
& 2016 & $ 89.7\% $ & $ 95.1\% $ & $ 88.6\% $ & $ 94.1\% $ & $ 87.7\% $ & $ 94.2\% $ \\ 
& 2017 & $ 90.0\% $ & $ 94.6\% $ & $ 87.8\% $ & $ 93.3\% $ & $ 86.6\% $ & $ 91.7\% $ \\ 
& 2018 & $ 90.4\% $ & $ 95.6\% $ & $ 87.5\% $ & $ 91.7\% $ & $ 89.0\% $ & $ 93.1\% $ \\
& average & $ 89.7\% $ & $ 94.8\% $ & $ 88.3\% $ & $ 93.4\% $ & $ 88.1\% $ & $ 93.0\% $ \\ 
\multicolumn{8}{c}{\multirow{2}{*}{quantile levels all $ = 50\% $}} \\ 
\\ 
& 2013 & $ 89.8\% $ & $ 94.2\% $ & $ 89.3\% $ & $ 95.2\% $ & $ 88.5\% $ & $ 93.3\% $ \\
& 2014 & $ 89.2\% $ & $ 95.3\% $ & $ 90.4\% $ & $ 95.4\% $ & $ 89.4\% $ & $ 94.8\% $ \\
& 2015 & $ 91.7\% $ & $ 96.0\% $ & $ 92.3\% $ & $ 95.7\% $ & $ 92.4\% $ & $ 95.9\% $ \\ 
& 2016 & $ 90.0\% $ & $ 95.5\% $ & $ 90.9\% $ & $ 95.5\% $ & $ 90.9\% $ & $ 94.9\% $ \\ 
& 2017 & $ 88.9\% $ & $ 95.2\% $ & $ 90.1\% $ & $ 96.0\% $ & $ 91.7\% $ & $ 95.9\% $ \\ 
& 2018 & $ 89.6\% $ & $ 94.9\% $ & $ 89.8\% $ & $ 95.4\% $ & $ 90.0\% $ & $ 95.4\% $ \\
& average & $ 89.9\% $ & $ 95.2\% $ & $ 90.5\% $ & $ 95.5\% $ & $ 90.5\% $ & $ 95.0\% $ \\ 
\bottomrule
\end{tabular}
%}
\end{table}

\begin{figure}[!hb] 
\centering
\caption{Confidence regions of the $ 20 $th percentiles of years 2013 and 2018 
by ELRT (solid), Wald (dashed), and nonparametric (dotted) methods, 
based on a simulated real data set of equal sample size $ n_{r} = 100 $. 
The true quantiles are marked with a diamond. The level of confidence is $ 95\% $.}
\label{RealData_CR_2013_2018_quan20_20}
\includegraphics[height = 241.5pt, width = 241.5pt]{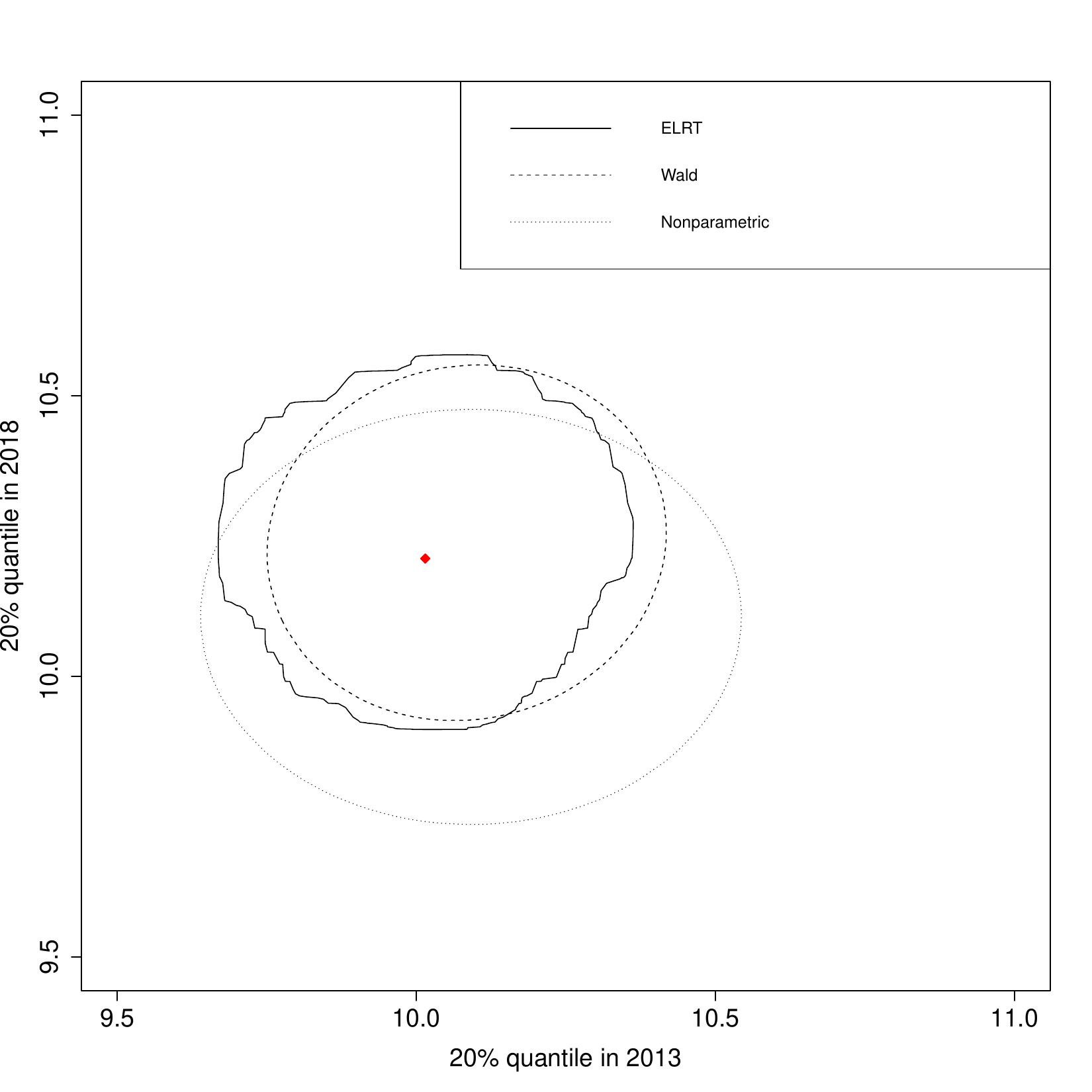} 
\end{figure}

In the next simulation, we focus on the confidence region of
the first quantiles of the household incomes in the years 2013 and 2018, 
namely the $ 20$th percentiles for these two years.
Figure~\ref {RealData_CR_2013_2018_quan20_20} shows the $ 95\% $ 
confidence regions using the three methods based on simulated real data of size $ n_{r} = 100 $. 
Table~\ref {RealData_Ava_Ecv2} gives the average coverages and areas of the three confidence regions, 
based on 1000 repetitions. 
The ELRT produces the most satisfactory confidence regions. 
The ELRT confidence regions improve the Wald confidence regions
by rightfully increased area to achieve more accurate coverage probabilities. 
They are much more efficient than the nonparametric confidence regions.

\begin{table}[!ht]
\centering
\caption{Empirical coverage probabilities and average areas for $ 20 $th 
percentiles of the years 2013 and 2018, based on real data of equal sample size.}
\label{RealData_Ava_Ecv2}
\resizebox{345.0pt}{!}{
\begin{tabular}{cccccc}
\toprule
\multirow{2}{*}{} & \multirow{2}{*}{Method} & \multicolumn{2}{c}{90\%} & \multicolumn{2}{c}{95\%}\\
\cmidrule(lr){3-4} \cmidrule(lr){5-6} 
& & Coverage probability & Area & Coverage probability & Area \\
\midrule
\multicolumn{6}{c}{\multirow{2}{*}{$ n_{r} = 100 $}} \\ 
\\ 
& ELRT & $ 89.00\% $ & $ 0.284 $ & $ 94.20\% $ & $ 0.379 $  \\ 
& Wald & $ 86.30\% $ & $ 0.245 $ & $ 91.80\% $ & $ 0.319 $  \\
& Nonparametric & $ 87.20\% $ & $ 0.358 $ & $ 91.60\% $ & $ 0.466 $  \\
\multicolumn{6}{c}{\multirow{2}{*}{$ n_{r} = 200 $}} \\ 
\\ 
& ELRT & $ 88.20\% $ & $ 0.130 $ & $ 93.40\% $ & $ 0.171 $  \\ 
& Wald & $ 86.10\% $ & $ 0.120 $ & $ 92.30\% $ & $ 0.156 $  \\
& Nonparametric & $ 88.80\% $ & $ 0.183 $ & $ 93.80\% $ & $ 0.238 $  \\
\bottomrule
\end{tabular}
}
\end{table}

\clearpage
\section{Power property and comparison}
\label{power}
Due to the linkage between the confidence region and the hypothesis test,
we are certain that the ELRT has superior power property based on
the simulation studies already done. 
At the same time, different tests have different higher power
regions in the space of the alternative hypotheses.
A generally inferior test can outperform other tests in specific regions.
We now use simulation to examine the power properties of the three tests.
We find the power properties do not vary much across different data types. 
To save space, we only present the simulation results based on real data. 

Consider the null hypothesis on values of the $ 20 $th percentiles of
years 2013 and 2018 with true values being $ \bxi_0 = (10.01, 10.21) $.
We examine the power of the three tests against
a range of false null hypotheses. One of them, for instance, is
\[
H_0: (\xi_1, \xi_2) = \bxi^* = \bxi_0 \times (0.99,  1.01) = (9.91, 10.31).
\]
We either inflate to deflate the true value by 1\% or 2\% leading to 8 false null hypotheses. 
We report the powers against these false $H_0$ 
in Table~\ref {RealData_Powers2} when the nominal
levels are 5\% and 10\% and the sample sizes are $n_r=100$ and $n_r=200$.

We observe that the rejection probabilities
of all three tests are above the corresponding nominal levels.
They increase when the sample size increases from $n_r=100$ to $n_r=200$.
These observations suggest the unbiasedness and consistency of the three tests.
We restrain from reading too much into some small differences as the sample size is
not sufficiently large. 
The power of ELRT is around 50\% when the assumed quantiles are 
2\% off from the truth and the sample size is $n_r = 100$ at level 5\%. 

Direct power comparison is most meaningful when tests under consideration
have the same size. 
Because we use the asymptotic distributions for all three methods,
there are non-ignorable differences in their null rejection rates
(see Table~\ref{RealData_Ava_Ecv2}). 
For each test, level, and sample size combination,
we calculate the average rejection rate.
The nonparametric test has lower power in general.
Yet the nonparametric test has higher rejection rates than ELRT 
2 out of 8 times when $n_r = 100$ at 5\% level.
However, the type I errors are 5.8\% and 8.4\%  
for ELRT and nonparametric in this case.
If this 44.8\% inflation factor in type I error is applied to their powers,
then ELRT would have higher powers in all 8 cases.
This general comment is applicable to all the other 3 combinations.

The Wald test seems to have higher power than the ELRT
on average in all 4 sample size and level combinations. 
However, its gain in lower type II error is at the cost of higher
type I error.  When $n_r=200$ and at level 10\%, 
Table~\ref{RealData_Ava_Ecv2} shows the ratio of their type I
errors is $1.178$.
Once we adjust the power of ELRT by this factor, the conclusion
will be reversed. 
This is the same for the other sample size and level combinations.

Although the three tests have different high power regions,
the ELRT is overall a better one.
The similar observations extend to unreported simulation results based on data generated from normal and gamma distributions.

\begin{table}[!ht]
\centering
\caption{Powers (in \%) for testing null hypotheses for
the $ 20 $th percentiles of years 2013 and 2018.}
\centerline{True value of the percentiles: $\bxi_0 = (10.01, 10.21)$.}
\label{RealData_Powers2}
\vspace{1ex}
\resizebox{345.0pt}{!}{
\begin{tabular}{cccccccc}
\toprule
&
& \multicolumn{2}{c}{ELRT} & \multicolumn{2}{c}{Wald}  & \multicolumn{2}{c}{Nonparametric} 
\\
\cmidrule(lr){3-4} \cmidrule(lr){5-6} \cmidrule(lr){7-8} 
\vspace{-2ex}
\\
 \multicolumn{2}{c}{Level of the test $\to$ } 
& $ 10\% $ & $ 5\% $ & $ 10\% $ & $ 5\% $ & $ 10\% $ & $ 5\% $ 
\\
\midrule
Change in scale $\downarrow$
&
$\bxi^*$ value in $H_0$ $\downarrow$
&
\multicolumn{6}{c}{\multirow{1}{*}{{\bf Rejection rates}}} \\ 
\midrule 
\multicolumn{8}{c}{\multirow{2}{*}{$ n_{r} = 100 $}} \\ 
\\ 
($-1, -1$)\%  & $ (9.91, 10.11) $ & $ 24.6 $ & $ 15.5 $ & $ 31.2 $ & $ 23.5 $ & $ 28.5 $ & $ 21.0 $ \\ 
($-1, +1$)\%  & $ (9.91, 10.31) $ & $ 23.1 $ & $ 14.0 $ & $ 25.0 $ & $ 15.2 $ & $ 20.6 $ & $ 13.4 $ \\ 
($+1, -1$)\%  & $ (10.12, 10.11) $ & $ 30.1 $ & $ 20.3 $ & $ 32.3 $ & $ 22.8 $ & $ 23.3 $ & $ 15.3 $ \\ 
($+1, +1$)\%  & $ (10.12, 10.31) $ & $ 26.6 $ & $ 15.0 $ & $ 18.0 $ & $ 9.3 $ & $ 12.3 $ & $ 6.4 $ \\ 
($-2, -2$)\%  & $ (9.81, 10.01) $ & $ 50.3 $ & $ 39.0 $ & $ 64.7 $ & $ 54.5 $ & $ 57.7 $ & $ 48.0 $ \\ 
($-2, +2$)\%  & $ (9.81, 10.41) $ & $ 51.7 $ & $ 40.7 $ & $ 62.3 $ & $ 49.5 $ & $ 46.5 $ & $ 34.7 $ \\ 
($+2, -2$)\%  & $ (10.22, 10.01) $ & $ 71.1 $ & $ 60.2 $ & $ 69.4 $ & $ 59.5 $ & $ 53.4 $ & $ 39.1 $ \\ 
($+2, +2$)\%  & $ (10.22, 10.41) $ & $ 62.9 $ & $ 52.5 $ & $ 56.9 $ & $ 41.4 $ & $ 37.4 $ & $ 22.9 $ \\ 
\multicolumn{2}{c}{average} & $ 42.6 $ & $ 32.2 $ & $ 45.0 $ & $ 34.5 $ & $ 35.0 $ & $ 25.1 $ \\ 
\multicolumn{8}{c}{\multirow{2}{*}{$ n_{r} = 200 $}} \\ 
\\ 
($-1, -1$)\%  & $ (9.91, 10.11) $ & $ 35.9 $ & $ 24.4 $ & $ 48.1 $ & $ 35.7 $ & $ 40.2 $ & $ 27.6 $ \\ 
($-1, +1$)\%  & $ (9.91, 10.31) $ & $ 34.8 $ & $ 23.3 $ & $ 41.6 $ & $ 27.9 $ & $ 28.0 $ & $ 18.6 $ \\ 
($+1, -1$)\%  & $ (10.12, 10.11) $ & $ 46.8 $ & $ 34.4 $ & $ 47.7 $ & $ 37.9 $ & $ 32.0 $ & $ 21.4 $ \\ 
($+1, +1$)\%  & $ (10.12, 10.31) $ & $ 40.7 $ & $ 28.5 $ & $ 33.9 $ & $ 20.9 $ & $ 22.9 $ & $ 13.4 $ \\ 
($-2, -2$)\%  & $ (9.81, 10.01) $ & $ 79.8 $ & $ 69.8 $ & $ 84.3 $ & $ 79.3 $ & $ 76.0 $ & $ 68.9 $ \\ 
($-2, +2$)\%  & $ (9.81, 10.41) $ & $ 80.1 $ & $ 71.0 $ & $ 85.1 $ & $ 79.7 $ & $ 73.2 $ & $ 62.3 $ \\ 
($+2, -2$)\%  & $ (10.22, 10.01) $ & $ 93.5 $ & $ 88.6 $ & $ 91.1 $ & $ 84.4 $ & $ 74.5 $ & $ 63.4 $ \\ 
($+2, +2$)\%  & $ (10.22, 10.41) $ & $ 90.4 $ & $ 82.5 $ & $ 85.4 $ & $ 76.8 $ & $ 69.9 $ & $ 56.3 $ \\ 
\multicolumn{2}{c}{average} & $ 62.8 $ & $ 52.8 $ & $ 64.7 $ & $ 55.3 $ & $ 52.1 $ & $ 41.5 $ \\ 
\bottomrule
\end{tabular}
}
\end{table}

%%%%%%%%%%%%%%%%%%%%%%%%%%%%%%%%%%%%%%%%%%%%%%%%%%%%%%%%%%%%%%%%%%%%%%

\newpage
\appendix
\section{Proofs of the main results} 
\label{App:Proofs}

This Appendix provides the proofs of the technical results. 
In the following proofs, without loss of generality, we proceed as if the sample proportions 
$ n_{k}/n $ do not depend on $ n $ and equal their limits $ \rho_{k} $. 
Our results are applicable as long as none of the populations have comparatively very small sample sizes.
Also, for the sake of convenience, with a generic function $ f (\by) $ we use 
\[ 
\frac{\partial f (\by^{*}) }{\partial \by} 
= 
\left. \frac{\partial f (\by) }{\partial \by} 
\right \rvert_{\by = \by^{*}}, ~~~
\frac{\partial^{2} f (\by^{*}) }{\partial \by \partial \by^{\top}} 
=
\left. \frac{\partial^{2} f (\by) }{\partial \by \partial \by^{\top}} 
\right \rvert_{\by = \by^{*}}. 
\] 
Moreover, the DRM parameters $ \btheta $ are arranged in the order 
\[
(
\theta_{1 1}, \theta_{2 1}, \ldots, \theta_{m 1}, 
\ldots, 
\theta_{1 2}, \theta_{2 2}, \ldots, \theta_{m 2}, 
\ldots, 
\theta_{1 d}, \theta_{2 d}, \ldots, \theta_{m d}
),  
\] 
where $ \theta_{i s} $ is the $ s $th component of the vector-valued parameter $ \btheta_{i} $. 
This order is needed for the expressions of the second derivative of $ \Dual (\blambda, \btheta) $ in the proof of Lemma \ref {lem:second_deriv} and for the covariance matrix of the first derivative in the proof of Lemma \ref {lem:first_deriv}.

%%%%%%%%%%%%%%%%%%%%%%%%%%%%%%%%%%%%%%%%%%%%%%%%%%%%%%%%%%%%%%%%%%%%%%

\subsection {Proof of Lemma \ref {lem:second_deriv}} 

This lemma asserts that the second derivative matrix of $ \Dual (\blambda, \btheta) $ has a finite and full-rank matrix as a limit.

\begin{proof} 

We first recognize that $ \Dual (\blambda, \btheta) $ can be written as a sum of $ m+1 $ sets of i.i.d. random variables: 
\begin{align} 
\Dual (\blambda, \btheta) = \sum_{k, j} \Dual_{k} (x_{k j}, \blambda, \btheta), 
 \label {dual_decomp}
\end{align} 
with 
\begin{align*}
\Dual_{k}  (x, \blambda, \btheta) 
= 
\btheta_{k}^{\top} \bq (x)
- \log \bar {h} (x, \btheta)
- \log \Big \{ 1 + \sum_{r \in I} \lambda_{r} \psi_{r} (x, \btheta) \Big \}. 
\end{align*}

Therefore, we may write 
\begin{align*} 
n^{-1} \frac { \partial^{2} \Dual (\bsm {0}, \btheta^{*}) } { \partial (\blambda, \btheta) \partial (\blambda, \btheta)^{\top} } 
= \sum_{k = 0}^{m} \rho_{k} \left [ n_{k}^{-1} \sum_{j = 1}^{n_{k}} \frac { \partial^{2} \Dual_{k} (x_{k j}, \bsm {0}, \btheta^{*}) } { \partial (\blambda, \btheta) \partial (\blambda, \btheta)^{\top} } \right ]. 
\end{align*} 
By the law of large numbers \citep {durrett2010probability}, as $ n \to \infty $, 
\begin{align*} 
n^{-1} 
\begin{pmatrix} 
\frac{ \partial^{2} \Dual (\bsm {0}, \btheta^{*}) }{\partial \blambda \partial \blambda^{\top}} 
& 
\frac{ \partial^{2} \Dual (\bsm {0}, \btheta^{*}) }{\partial \blambda \partial \btheta^{\top}} \\ 
~ \\ 
\frac{ \partial^{2} \Dual (\bsm {0}, \btheta^{*}) }{\partial \btheta \partial \lambda^{\top}} 
&
\frac{ \partial^{2} \Dual (\bsm {0}, \btheta^{*}) }{\partial \btheta \partial \btheta^{\top}} 
\end{pmatrix}
\to
\begin{pmatrix} 
S_{\blambda \blambda} & S_{\blambda \btheta} \\ 
S_{\btheta \blambda} & S_{\btheta \btheta} 
\end{pmatrix}, 
\end{align*} 
for some block matrix $ S $ given by 
\begin{align*} 
S = \sum_{k = 0}^{m} \rho_{k} \mathbbm {E}_{k} \left [ \frac { \partial^{2} \Dual_{k} (X, \bsm {0}, \btheta^{*}) } { \partial (\blambda, \btheta) \partial (\blambda, \btheta)^{\top} } \right ]. 
\end{align*}  
Here we remark again that we assume that the sample proportions $ n_{k}/n $ do not change with $ n $ and always equal their limits $ \rho_{k} $. 

Next, we show that $ S $ has full rank. 
We first give the following expressions:
\begin{align*} 
\dfrac { \partial^{2} \Dual_{k} (x, \bsm {0}, \btheta^{*}) } 
{ \partial \blambda \partial \blambda^{\top} } 
& = 
\bpsi (x, \btheta^{*}) \bpsi^{\top} (x, \btheta^{*}), 
\\ 
\dfrac { \partial^{2}  \Dual_{k} (x, \bsm {0}, \btheta^{*}) } 
	{ \partial \btheta \partial \btheta^{\top} } 
& = 
[ \bq (x) \bq^{\top} (x) ] \otimes \left [ \bh (x, \btheta^{*}) \bh^{\top} (x, \btheta^{*}) - \mathrm {diag} \{ \bh (x, \btheta^{*}) \} \right ], 
\\ 
\frac { \partial^{2} \Dual_{k} (x, \bsm {0}, \btheta^{*}) } 
	{ \partial \blambda \partial \btheta^{\top} } 
& = 
 \bq^{\top} (x) \otimes 
 \left [
 \bpsi (x, \btheta^{*}) \bh^{\top} (x, \btheta^{*}) - 
 \mathrm {diag} \{ \bpsi (x, \btheta^{*}) \} 
\begin{pmatrix} 
\be_{I_{1}} & \cdots & \be_{I_{l}}
\end{pmatrix}^{\top}
 \right ], 
\end{align*} 
where $ \otimes $ is the Kronecker product, $ \be_{i} $ is a vector of length $ m $ that has 1 in the $ i $th entry and 0 elsewhere (we define $ \be_{0} = \bsm {0} $ by convention), and $ I_{j} $ is the population index of the $ j $th quantile of interest. 

Based on the above expressions, we first note that
\begin{align*} 
S_{\btheta \btheta} = 
- \sum_{k = 0}^{m} \rho_{k} \bbE_{k} 
\left [ 
\{ \bq (X) \otimes [\be_{k} - \bh (X, \btheta^{*})] \} 
         \{ \bq (X) \otimes [\be_{k} - \bh (X, \btheta^{*})] \}^{\top} 
\right ], 
\end{align*} 
which is clearly negative semidefinite.
We now strengthen the conclusion to negative definite.
By Condition \eqref {Condition.ii},  $ \bbE_{0} [ \bq (X) \bq^{\top} (X) ] $ is positive definite. 
Since 
$ h_{r} (x, \btheta^{*}) = \exp (\btheta_{r}^{\top} \bq (x))$, we have that 
\[
\bbE_{k} 
\left [ 
\{ \be_{k} - \bh (X, \btheta^{*}) \} \{ \be_{k} - \bh (X, \btheta^{*}) \}^{\top} 
\right ] 
\]
is positive definite.
Simple algebra leads to the negative definiteness of $ S_{\btheta \btheta} $.
For the same reason, $ S_{\blambda \blambda} $ is positive definite if $ \bpsi (x, \btheta^{*}) $ does not degenerate, which is assured because
\[
\varphi_{r} (x, \btheta, \bxi) 
= 
h_{r} (x, \btheta) [ \mathbbm {1} (x \leq \xi_{r}) - \tau_{r} ].
\]

From
\[
\begin{pmatrix} 
\mathbbm {I} & - S_{\blambda \btheta}S^{-1}_{\btheta \btheta} \\ 
\bsm {0} & \mathbbm {I}
\end{pmatrix}
\times
\begin{pmatrix} 
S_{\blambda \blambda} & S_{\blambda \btheta} \\ 
S_{\btheta \blambda} & S_{\btheta \btheta} 
\end{pmatrix}
=
\begin{pmatrix} 
S_{\blambda \blambda} - S_{\blambda \btheta} S_{\btheta \btheta}^{-1} S_{\btheta \blambda } & \bsm {0} \\ 
S_{\btheta \blambda}  & S_{\btheta \btheta} 
\end{pmatrix}, 
\]
we conclude that $ S $ has full rank if  
$ S_{\blambda \blambda} - S_{\blambda \btheta} S_{\btheta \btheta}^{-1} S_{\blambda \btheta}^{\top} $ does. 
Because $ S_{\blambda \blambda} $ is positive definite and $ S_{\btheta \btheta}^{-1} $ is negative definite, $ S_{\blambda \blambda} - S_{\blambda \btheta} S_{\btheta \btheta}^{-1} S_{\blambda \btheta}^{\top} $ must be positive definite, and so it has full rank. 
This completes the proof that $ S $ has full rank. 

\end{proof}

%%%%%%%%%%%%%%%%%%%%%%%%%%%%%%%%%%%%%%%%%%%%%%%%%%%%%%%%%%%%%%%%%%%%%%

\subsection {Proof of Lemma \ref {lem:first_deriv}}

\begin{proof} 

The first conclusion of this lemma is that the first derivative of $ \Dual (\blambda, \btheta) $ in \eqref{new_dual1} has zero expectation when evaluated at $ (\bsm {0}, \btheta^{*}) $. 
Recall that
\[
\Dual (\blambda, \btheta) 
=
 \sum_{k, j} \btheta_{k}^{\top} \bq (x_{k j})
 -
\sum_{k, j}  \log \bar {h} (x_{k j}, \btheta) 
- \sum_{k, j} 
\log  
\big \{ 1 + \sum_{r \in I} \lambda_{r} \psi_{r} (x_{k j}, \btheta) \big \}.
\]
For any $ r \in I $, the partial derivative of $ \Dual (\blambda, \btheta) $ with respect to $ \blambda_{r} $ is given by 
\[
\frac { \partial \Dual (\blambda, \btheta) } { \partial \lambda_{r} } 
= 
- \sum_{k, j} 
\frac { \psi_{r} (x_{k j}, \btheta) } 
{1 + \sum_{i \in I} \lambda_{i} \psi_{i} (x_{k j}, \btheta) }. 
\]
At $ \blambda^{*} = \bsm {0} $ and $ \btheta = \btheta^{*} $, this reduces to
\[
\frac { \partial \Dual (\bsm {0}, \btheta^{*})  } { \partial \lambda_{r} } 
= 
- \sum_{k, j} \psi_{r} (x_{k j}, \btheta^{*}).  
\]
Hence, we have 
\begin{align} 
\mathbbm {E} 
\left [ \frac { \partial \Dual (\bsm {0}, \btheta^{*})  } { \partial \lambda_{r} } \right ]  
& = 
- \sum_{k, j} \int \psi_{r} (x, \btheta^{*}) \mathrm {d} G_{k} (x) 
\nonumber \\ 
& = 
- \int \psi_{r} (x, \btheta^{*}) 
\big \{ \sum_{k = 0}^{m} n_{k} h_{k} (x, \btheta^{*}) \big \} \mathrm {d} G_{0} (x)
\nonumber \\ 
& = - n \int \varphi_{r} (x, \btheta^{*}, \bxi^{*}) \mathrm {d} G_{0} (x) 
= 0. 
\label{App.eq1.1}
\end{align} 

For each $ i = 1, 2, \ldots, m, \, s = 1, 2, ... , d $, and 
at $ \blambda = \bsm {0} $ and $ \btheta = \btheta^{*} $, we have 
\[
\frac { \partial \Dual (\bsm {0}, \btheta^{*}) } { \partial \theta_{i s} } 
 = 
 \sum_{j = 1}^{n_{i}} q_{s} (x_{i j}) 
 -  
 \sum_{k, j} \rho_{i} q_{s} (x_{k j}) h_{i} (x_{k j}, \btheta^{*})/\bar{h} (x_{k j}, \btheta^{*}).
\]
For the first term, it can be seen that
\begin{align*} 
\mathbbm {E} \Big [ \sum_{j = 1}^{n_{i}} q_{s} (x_{i j})  \Big ]
=
n_{i} \int q_{s} (x) h_{i} (x, \btheta^{*}) \mathrm {d} G_{0} (x) .
\end{align*} 
At the same time, for the second term, we have
\begin{align*} 
\mathbbm {E} 
\Big [  
\sum_{k, j} \rho_{i} q_{s} (x_{k j}) h_{i} (x_{k j}, \btheta^{*})/\bar{h} (x_{k j}, \btheta^{*})
\Big ] 
& = 
n_{i} \int q_{s} (x) h_{i} (x, \btheta^{*}) \mathrm {d} G_{0} (x).
\end{align*} 
Therefore, we find that
\begin{align} 
\mathbbm {E} \Big [
\frac { \partial \Dual (\bsm {0}, \btheta^{*}) } { \partial \theta_{i s} } 
\Big ] = 0. 
\label{App.eq1.2}
\end{align} 

Combining \eqref{App.eq1.1} and \eqref{App.eq1.2}, we conclude that
\begin{align*} 
\mathbbm {E} \Big [
\frac { \partial \Dual (\bsm {0}, \btheta^{*}) } { \partial (\blambda, \btheta) } 
\Big ] = \bsm {0}. 
\end{align*}

The second conclusion of this lemma is the asymptotic normality of the first derivative. 
Despite its complex expression, we can see that $ \partial \Dual (\blambda, \btheta)/\partial (\blambda, \btheta) $ is a sum of $ m+1 $ sets of i.i.d. random variables of sizes $ n_{r} = n \rho_{r} $ with mean zero and finite second moment in the matrix sense. 
Recall \eqref {dual_decomp} from the proof of Lemma \ref {lem:second_deriv} that 
\begin{align*} 
\Dual (\blambda, \btheta) = \sum_{k, j} \Dual_{k} (x_{k j}, \blambda, \btheta), 
\end{align*} 
where 
\begin{align*}
\Dual_{k}  (x, \blambda, \btheta) 
= 
\btheta_{k}^{\top} \bq (x)
- \log \bar {h} (x, \btheta)
- \log \Big \{ 1 + \sum_{r \in I} \lambda_{r} \psi_{r} (x, \btheta) \Big \}. 
\end{align*}

We may write
\begin{align*}
\frac {\partial \Dual (\blambda, \btheta)}{\partial (\blambda, \btheta)}
= 
\sum_{k = 0}^{m} 
\left \{  
\sum_{j = 1}^{n_{k}} 
\left [ \frac {\partial \Dual_{k} (x_{k j}, \blambda, \btheta)}{\partial (\blambda, \btheta)} 
-
\mathbbm {E}_{k} 
\left (
\frac {\partial \Dual_{k} (X, \blambda, \btheta)}{\partial (\blambda, \btheta)}
\right ) \right ] 
\right \}.  
\end{align*}

For each $ k = 0, 1, \ldots, m $, as $ n_{k} \to \infty $, 
\begin{align*} 
T_{k} \coloneqq n_{k}^{-1/2} \sum_{j = 1}^{n_{k}} \left [ \frac { \partial \Dual_{k} (x_{k j}, \bsm {0}, \btheta^{*}) } { \partial (\blambda, \btheta) } - \mathbbm {E}_{k} \left (\frac {\partial \Dual_{k} (X, \bsm {0}, \btheta^{*})}{\partial (\blambda, \btheta)} \right ) \right ] 
\end{align*} 
has a limiting distribution of normal with mean zero and finite second moment in the matrix sense, by the multivariate central limit theorem for triangular arrays \citep {durrett2010probability}. 
Because $ T_{0}, T_{1}, \ldots, T_{m} $ are independent of each other, the targeted quantity  
\begin{align*}
n^{-1/2} \frac {\partial \Dual (\bsm {0}, \btheta^{*})}{\partial (\blambda, \btheta)}
= 
\sum_{k = 0}^{m} \rho_{k}^{1/2} T_{k} 
\end{align*} 
is asymptotically normal with mean zero. 

We now give the expression $ V $ for the covariance matrix in the limiting distribution.
Let $ V_{k} $ be the asymptotic covariance matrix of $ T_{k} $, then we have 
\begin{align*} 
V = \sum_{k = 0}^{m} \rho_{k} V_{k}.  
\end{align*} 
The expression for $ V_{k} $ is given by 
\begin{align*} 
V_{k} 
= & 
\mathbbm {E}_{k} \left [ \left ( \frac { \partial \Dual_{k} (X, \bsm {0}, \btheta^{*}) } { \partial (\blambda, \btheta) } \right ) \left ( \frac { \partial \Dual_{k} (X, \bsm {0}, \btheta^{*}) } { \partial (\blambda, \btheta) } \right )^{\top} \right ] \\ 
& - \mathbbm {E}_{k} \left [ \frac {\partial \Dual_{k} (X, \bsm {0}, \btheta^{*})}{\partial (\blambda, \btheta)} \right ] \mathbbm {E}_{k} \left [ \frac {\partial \Dual_{k} (X, \bsm {0}, \btheta^{*})}{\partial (\blambda, \btheta)} \right ]^{\top}. 
\end{align*}  
After some algebra, we find that 
\begin{gather*} 
\frac {\partial \Dual_{k} (x, \bsm {0}, \btheta^{*})} {\partial \blambda} = - \bpsi (x, \btheta^{*}), \\ 
\frac {\partial \Dual_{k} (x, \bsm {0}, \btheta^{*})} {\partial \btheta} = \bq (x) \otimes [ \be_{k} - \bh (x, \btheta^{*}) ],  
\end{gather*}   
where $ \be_{k} $ is a unit vector with the $ k $th element being 1 ($ \be_{0} = \bsm {0} $ by convention).  
We have 
\begin{align*} 
\sum_{k = 0}^{m} \rho_{k} \mathbbm {E}_{k} \left [ \left ( \frac { \partial \Dual_{k} (X, \bsm {0}, \btheta^{*}) } { \partial (\blambda, \btheta) } \right ) \left ( \frac { \partial \Dual_{k} (X, \bsm {0}, \btheta^{*}) } { \partial (\blambda, \btheta) } \right )^{\top} \right ] = 
\begin{pmatrix} 
S_{ \blambda \blambda } & \bsm {0} \\ 
\bsm {0} & - S_{ \btheta \btheta } \\ 
\end{pmatrix}. 
\end{align*} 
Let 
\[ 
W = 
\begin{pmatrix}  
\rho_{0}^{-1} \bone_{m} \bone_{m}^{\top} + \mathrm {diag} \{ \rho_{1}^{-1}, \ldots, \rho_{m}^{- 1} \} & \bsm {0} \\ 
\bsm {0} & \bsm {0} \\ 
\end{pmatrix}, 
\] 
where $ \bone_{m} $ is an $ m $-dimensional vector of ones; 
we then also have 
\begin{align*} 
\sum_{k = 0}^{m} \rho_{k} \mathbbm {E}_{k} \left [ \frac {\partial \Dual_{k} (X, \bsm {0}, \btheta^{*})}{\partial (\blambda, \btheta)} \right ] \mathbbm {E}_{k} \left [ \frac {\partial \Dual_{k} (X, \bsm {0}, \btheta^{*})}{\partial (\blambda, \btheta)} \right ]^{\top} = 
S 
\begin{pmatrix} 
\bsm {0} & \bsm {0} \\ 
\bsm {0} & W \\
\end{pmatrix}
S. 
\end{align*} 

Finally, we get 
\begin{align*} 
V = 
\begin{pmatrix} 
S_{ \blambda \blambda } & \bsm {0} \\ 
\bsm {0} & - S_{ \btheta \btheta } \\ 
\end{pmatrix} 
- 
S 
\begin{pmatrix} 
\bsm {0} & \bsm {0} \\ 
\bsm {0} & W \\
\end{pmatrix}
S. 
\end{align*} 
This completes the proof that $ n^{-1/2} \partial \Dual (\bsm {0}, \btheta^{*})/\partial (\blambda, \btheta) $ is asymptotically normal. 

\end{proof}

%%%%%%%%%%%%%%%%%%%%%%%%%%%%%%%%%%%%%%%%%%%%%%%%%%%%%%%%%%%%%%%%%%%%%%

\subsection {Proof of Lemma \ref {lem:para_normality}}

\begin{proof} 

Given $ \btheta $, let $ \blambda (\btheta) $ be the solution to 
\begin{align*} 
\sum_{k, j} \frac {\bpsi (x_{k j}, \btheta)} {1 +  \blambda^{\top} \bpsi (x_{k j}, \btheta)} = \bsm {0}. 
\end{align*} 
We first prove that uniformly for any $ \btheta $ in the $ n^{-1/3} $-neighbourhood of $ \btheta^{*} $, 
$ \blambda (\btheta) $ is $ O (n^{-1/3}) $. 
For notational convenience, in this section we omit $ \btheta $ in $ \blambda (\btheta) $ if this does not cause any confusion. 

Following the typical proof in \citet {owen2001empirical}, 
the claim is true if uniformly for $ \btheta $ such that $ \| \btheta - \btheta^{*} \| \leq n^{-1/3} $ we have 
\begin{enumerate}[(i)]

\item
$ \sum_{k, j} \bpsi (x_{k j}, \btheta) = O (n^{2/3}) $; 
\label {NBDtarget1} 

\item
$ n^{-1} \sum_{k, j} \bpsi (x_{k j}, \btheta) \bpsi^{\top} (x_{k j}, \btheta) $ has a positive definite limit. 
\label {NBDtarget2} 

\end{enumerate}
We omit other details but prove the above results. 
Note that little $ o $ and big $ O $ without $ p $ in the subscript are orders in the sense of almost surely.

Recall that $ \sum_{k} \mathbbm {E}_{k} [\bpsi (X, \btheta^{*})] = \bsm {0} $ 
as shown in Lemma \ref {lem:first_deriv}.
We have 
\begin{align} 
\label{eq20b}
\sum_{k, j} \bpsi (x_{k j}, \btheta^{*}) 
& = \sum_{k} \Big \{ \sum_{j} \bpsi (x_{k j}, \btheta^{*}) - \mathbbm {E}_{k} [\bpsi (X, \btheta^{*})] \Big \} 
\nonumber \\ 
& =  \sum_{k}  \{ O (\sqrt {n \log \log n}) \} = O (n^{2/3}),
\end{align} 
applying the law of the iterated logarithm to each $ k $.

For $ \btheta $ in a small neighbourhood of $ \btheta^{*} $, 
there is a generic nonrandom constant $ C $ such that
\begin{align} 
\label{eq20c}
\sum_{k, j} \| \partial \bpsi (x_{k j}, \btheta) / \partial \btheta \|
\leq 
C \sum_{k, j} \| \bq (x_{k j}) \| = O (n), 
\end{align}  
with the order in the last step derived from the finite moment assumption on $ \bq (X) $.
Applying \eqref {eq20b} and \eqref {eq20c}, with $ \bar \btheta $ being a value between $ \btheta $ and $ \btheta^{*} $,
we get
\begin{align*} 
\sum_{k, j} \bpsi (x_{k j}, \btheta)
& = \sum_{k, j} \bpsi (x_{k j}, \btheta^{*}) 
+
\sum_{k, j} \frac {\partial \bpsi (x_{k j}, \bar \btheta)} {\partial \btheta} (\btheta - \btheta^{*})
\nonumber \\
& = 
\sum_{k, j} \frac {\partial \bpsi (x_{k j}, \bar \btheta)} {\partial \btheta} (\btheta - \btheta^{*}) 
+ O (n^{2/3}) = O (n^{2/3}). 
\end{align*} 
This proves \eqref {NBDtarget1}.

Recall that $ \bpsi (x, \btheta) = \{ \psi_{r} (x, \btheta) : r \in I \} $ and observe
\[ 
| \psi_r (x, \btheta) |
= \left | \frac {h_{r} (x, \btheta)}{\sum_{r = 0}^{m} \rho_{r} h_{r} (x, \btheta)} [ \mathbbm {1} (x \leq \xi_{r}) - \tau_{r} ] \right | \leq \rho_r^{-1} = O(1).
\]
By focusing on $ \btheta $ in an $ n^{-1/3} $-neighbourhood of $ \btheta^{*} $,
we have
\begin{align*} 
&
\sum_{j = 1}^{n_{k}} \{ \bpsi (x_{k j}, \btheta) \bpsi^{\top} (x_{k j}, \btheta)
-
\bpsi (x_{k j}, \btheta^{*}) \bpsi^{\top} (x_{k j}, \btheta^{*}) \} \\
&=
\sum_{j = 1}^{n_{k}} 
\Big [ \bpsi (x_{k j}, \btheta^{*}) \{ \bpsi (x_{k j}, \btheta)- \bpsi (x_{k j}, \btheta^{*}) \}^{\top}
+
\{ \bpsi (x_{k j}, \btheta)- \bpsi (x_{k j}, \btheta^{*}) \}  \bpsi^{\top} (x_{k j}, \btheta) \Big ] 
\\
& \leq
\{ \max_{k, j} \sup_{\btheta} \bpsi (x_{k j}, \btheta) \}
\sum_{j = 1}^{n_{k}} \| \bpsi (x_{k j}, \btheta) - \bpsi (x_{k j}, \btheta^{*}) \|  \\
&= O (n^{2/3}) = o (n). 
\end{align*} 

Therefore, we have
\begin{align*} 
n^{-1} \sum_{k, j} \{ \bpsi (x_{k j}, \btheta) \bpsi^{\top} (x_{k j}, \btheta) \}
& =
n^{-1} \sum_{k, j} \{ \bpsi (x_{k j}, \btheta^{*}) \bpsi^{\top} (x_{k j}, \btheta^{*}) \}+ o (1) \\ 
& \to S_{\blambda \blambda}, 
\end{align*} 
which is clearly positive definite. 
This proves \eqref {NBDtarget2}.

As we have remarked, the validity of \eqref {NBDtarget1} and \eqref {NBDtarget2} implies that uniformly for $ \btheta - \btheta^{*} = O (n^{-1/3}) $, 
\begin{align} 
\blambda (\btheta) = O (n^{-1/3}).
\label {lambda_order}
\end{align}
Following the same line of the proof, we also have a stronger order for $ \lambda (\btheta) $ when $ \btheta = \btheta^{*} $: 
\begin{align} 
\blambda (\btheta^{*}) = o (n^{-1/3}). 
\label {lambda_order_true_theta}
\end{align} 

The next stage of the proof is dedicated to showing that $ \hat \btheta - \btheta^{*} = O (n^{-1/3}) $. 
We consider a function of $ \btheta $:  
\begin{align*} 
L (\btheta) = \Dual (\blambda (\btheta), \btheta). 
\end{align*} 
It can easily be seen that $ \hat \btheta $ is a maximizer of $ L (\btheta) $. 
Since $ L (\btheta) $ is a smooth function, there must be a maximizer of $ L (\btheta) $ in the compact set $ \{ \btheta: \| \btheta - \btheta^{*} \| \leq n^{-1/3} \} $. 
We prove that this maximizer is attained in the interior of the compact set by showing that $ L (\btheta) < L (\btheta^{*}) $ uniformly for $ \btheta $ on the boundary of the compact set. 
For any unit vector $ \ba $ and $ \btheta = \btheta^{*} + n^{-1/3} \ba $, expanding $ L (\btheta) $ at $ \btheta^{*} $ yields (see \citet {folland2002advanced})  
\begin{align} 
L (\btheta) = L (\btheta^{*}) 
+ n^{-1/3} \frac {\partial L (\btheta^{*})} {\partial \btheta} \ba 
+ n^{-2/3} \ba^{\top} \frac {\partial^{2} L (\btheta^{*})} {\partial \btheta \partial \btheta^{\top}} \ba + \varepsilon_{n},  
\label {L_Expansion_1} 
\end{align} 
where $ \varepsilon_{n} $ is the Lagrange remainder term in obvious notation: 
\begin{align*} 
\varepsilon_{n} = \frac {1} {6} n^{-1} \sum_{ |\alpha| = 3 } \partial^{\alpha} L (\underline {\btheta}) \, \, \ba^{\alpha}, 
\end{align*} 
for some $ \underline {\btheta} $ between $ \btheta^{*} $ and $ \btheta $. 
By the uniform boundedness of the third-order derivatives of $ L (\btheta) $, we have $ \varepsilon_{n} = O (1) $ uniformly over $ \ba $. 

For the first term in the expansion, we note that $ \blambda (\btheta^{*}) = o (n^{-1/3}) $ as given in \eqref {lambda_order_true_theta}, and this implies
\[
\frac {\partial \Dual ( \blambda (\btheta^{*}), \btheta^{*})} {\partial \btheta} 
= \frac {\partial \Dual (\bsm {0}, \btheta^{*})} {\partial \btheta} + O (n) (\blambda (\btheta^{*}) - \bsm {0})
= o (n^{2/3}), 
\]
with the order of $ \partial \Dual (\bsm {0}, \btheta^{*})/{\partial \btheta} $ implied by Lemma \ref {lem:first_deriv}.
Therefore,
\begin{align*} 
\frac {\partial L (\btheta^{*})} {\partial \btheta} 
& = 
\frac {\partial \Dual (\blambda (\btheta^{*}), \btheta^{*})} 
{\partial \blambda} \frac {\partial \blambda (\btheta^{*})} {\partial \btheta} 
+ \frac {\partial \Dual (\blambda (\btheta^{*}), \btheta^{*})} {\partial \btheta} \\
= &
\bsm {0} + \frac {\partial \Dual (\blambda (\btheta^{*}), \btheta^{*})} {\partial \btheta}
= o (n^{2/3}).
\end{align*}

For the second term in the expansion, we proceed as follows.
With $ \blambda (\btheta^{*}) = o (n^{-1/3}) $ as given in \eqref {lambda_order_true_theta} and Lemma \ref {lem:second_deriv},
we first note that
\[
\frac {\partial^{2} \Dual ( \blambda (\btheta^{*}), \btheta^{*})} 
{\partial (\blambda, \btheta) \partial (\blambda, \btheta)^{\top}} 
= 
\frac {\partial^{2} \Dual (\bsm {0}, \btheta^{*})} 
{\partial (\blambda, \btheta) \partial (\blambda, \btheta)^{\top}} + o (n^{2/3}) 
= n [S + o (1)].  
\]
Taking derivatives with respect to $ \btheta $ on both sides of the identity
\[
\frac{\partial \Dual (\blambda (\btheta), \btheta) }{ \partial \blambda} = \bsm {0}, 
\]
and then setting $ \btheta = \btheta^{*} $, we further have 
\begin{align*} 
\frac {\partial  \blambda (\btheta^{*})} {\partial \btheta} 
= 
- \left [ \frac {\partial^{2} \Dual ( \blambda (\btheta^{*}), \btheta^{*})} 
{\partial \blambda \partial \blambda^{\top}} \right ]^{-1} 
\left [ \frac {\partial^{2} \Dual ( \blambda (\btheta^{*}), \btheta^{*})} 
{\partial \blambda \partial \btheta^{\top}} \right ]. 
\end{align*} 
Hence,
\begin{align*} 
\frac {\partial^{2} L (\btheta^{*})} {\partial \btheta \partial \btheta^{\top}} 
= &
\frac {\partial^{2} \Dual (\blambda (\btheta^{*}), \btheta^{*})} 
{\partial \btheta \partial \blambda^{\top}} \frac {\partial \blambda (\btheta^{*})} {\partial \btheta} 
+ 
\frac {\partial^{2} \Dual (\blambda (\btheta^{*}), \btheta^{*})} {\partial \btheta \partial \btheta^{\top}}
\\
= & 
- \left [ \frac {\partial^{2} \Dual ( \blambda (\btheta^{*}), \btheta^{*})} {\partial \blambda \partial \btheta^{\top}} \right ]^{\top} 
\left [ \frac {\partial^{2} \Dual ( \blambda (\btheta^{*}), \btheta^{*})} {\partial \blambda \partial \blambda^{\top}} \right ]^{-1} 
\left [ \frac {\partial^{2} \Dual ( \blambda (\btheta^{*}), \btheta^{*})} {\partial \blambda \partial \btheta^{\top}} \right ] \\  
& + 
\frac {\partial^{2} \Dual ( \blambda (\btheta^{*}), \btheta^{*})} {\partial \btheta \partial \btheta^{\top}}
\\
= &
n [ -S_{\blambda \btheta}^{\top} S_{\blambda \blambda}^{-1} S_{\blambda \btheta} + S_{\btheta \btheta} + o (1) ]. 
\end{align*} 

Therefore, the expansion of $ L (\btheta) $ in \eqref {L_Expansion_1} becomes
\begin{align*} 
L (\btheta) - L (\btheta^{*}) 
= 
n^{1/3} \ba^{\top} 
\{ - S_{\blambda \btheta}^{\top} S_{\blambda \blambda}^{-1} S_{\blambda \btheta}
+ S_{\btheta \btheta}
\}
\ba + o (n^{1/3}). 
\end{align*} 
The matrix in the quadratic form is negative definite, following the line of an argument in the proof of Lemma \ref {lem:second_deriv}.
Hence, as $ n \to \infty $, with probability 1, 
\begin{align*} 
L (\btheta^{*} + n^{-1/3} \ba) < L (\btheta^{*}), 
\end{align*} 
uniformly over all unit vector $ \ba $. 
This proves  
\begin{align*} 
\hat \btheta - \btheta^{*} = O (n^{-1/3}), 
\end{align*} 
and together with \eqref {lambda_order} further implies that
\begin{align*} 
\hat \blambda = O (n^{-1/3}). 
\end{align*}

%======================================================================= 

We are now ready to prove the asymptotic normality of $ (\hat \blambda, \hat \btheta) $. 
Expanding $ \partial \Dual (\hat \blambda, \hat \btheta) / \partial (\blambda, \btheta) $ at $ (\bsm {0}, \btheta^{*}) $, we get  
\begin{align*} 
\bsm {0} = \frac { \partial \Dual (\hat \blambda, \hat \btheta) } { \partial (\blambda, \btheta) } 
& = \frac { \partial \Dual (\bsm {0}, \btheta^{*}) } { \partial (\blambda, \btheta) } + \frac { \partial^{2} \Dual (\bsm {0}, \btheta^{*}) } { \partial (\blambda, \btheta) \partial (\blambda, \btheta)^{\top} } 
\begin{pmatrix}
\hat \blambda - \bsm {0} \\ 
\hat \btheta - \btheta^{*} 
\end{pmatrix}
+ O (n^{1/3}). 
%\label {Taylor_E3}
\end{align*}  

By Lemmas \ref {lem:second_deriv} and \ref {lem:first_deriv}, we get
\begin{align} 
\sqrt {n} 
\begin{pmatrix}
\hat \blambda - \bsm {0} \\
\hat \btheta - \btheta^{*}
\end{pmatrix}
& = - S^{-1} \left [ n^{-1/2} 
\frac { \partial \Dual (\bsm {0}, \btheta^{*}) } { \partial (\blambda, \btheta) } \right ] 
+ o_p (1)
\overset {d} \to N (\bsm {0}, S^{-1} V S^{-1}), 
\label{eq20e}
\end{align} 
as $ n \to \infty $. 
 
\end{proof}

%%%%%%%%%%%%%%%%%%%%%%%%%%%%%%%%%%%%%%%%%%%%%%%%%%%%%%%%%%%%%%%%%%%%%%

\subsection{Proof of Theorem \ref {thm:ELRT_Stat_chisq}} 

\begin{proof} 

We notice that, as shown in \citet {cai2017hypothesis}, 
\[ 
\sup_{ \btheta, G_{0} } \{ \ell_{n} (\btheta, G_{0}) \} = \sup_{\btheta} \Dual (\bsm {0}, \btheta)  - n \log n.  
\] 
From \eqref {dual_relation} we also have 
\begin{align*} 
\tilde \ell_{n} (\bxi^{*}) = \Dual (\hat \blambda, \hat \btheta) - n \log n. 
\end{align*} 
These relations lead to
\begin{align} 
R_{n} 
& = 
2 \left [ \sup_{\btheta} \Dual (\bsm {0}, \btheta) - \Dual (\hat \blambda, \hat \btheta) \right ] \nonumber \\
&=
2 \left [ \sup_{\btheta} \Dual (\bsm {0}, \btheta) - \Dual (\bsm {0}, \btheta^{*}) \right ] 
- 
2 \left [ \Dual (\hat \blambda, \hat \btheta) - \Dual (\bsm {0}, \btheta^{*}) \right ].
\label {ELRT_decompose}
\end{align} 

\citet {cai2017hypothesis} show in the proof of their Theorem 1 that 
\begin{align*} 
\sup_{\btheta} \Dual (\bsm {0}, \btheta ) - \Dual (\bsm {0}, \btheta^{*}) 
& = - \frac {1} {2} 
\left [ n^{-1/2} \frac { \partial \Dual (\bsm {0}, \btheta^{*}) } { \partial \btheta } \right ]^{\top} 
S_{\btheta \btheta}^{-1} 
\left [ n^{-1/2} \frac { \partial \Dual (\bsm {0}, \btheta^{*}) } { \partial \btheta } \right ] 
+ o_{p} (1),  
\end{align*} 
for the same $ S_{\btheta \btheta} $ given in the proof of Lemma \ref {lem:second_deriv}.  

For the second term in \eqref {ELRT_decompose}, utilizing the expansion of $ \hat \blambda $ and $ \hat \btheta - \btheta^{*} $ given in 
\eqref {eq20e}, we have
\begin{align*}  
& \Dual (\hat \blambda, \hat \btheta) - \Dual (\bsm {0}, \btheta^{*}) \\ 
= &
\frac { \partial \Dual (\bsm {0}, \btheta^{*}) } { \partial (\blambda, \btheta) } 
\begin{pmatrix}
\hat \blambda - \bsm {0} \\ 
\hat \btheta - \btheta^{*} 
\end{pmatrix} 
+ 
\frac {1} {2}
\begin{pmatrix}
\hat \blambda - \bsm {0} \\ 
\hat \btheta - \btheta^{*} 
\end{pmatrix}^{\top} 
\frac { \partial^{2} \Dual (\bsm {0}, \btheta^{*}) } 
{ \partial (\blambda, \btheta) \partial (\blambda, \btheta)^{\top} } 
\begin{pmatrix}
\hat \blambda - \bsm {0} \\ 
\hat \btheta - \btheta^{*} 
\end{pmatrix}
+ o_{p} (1)  \\ 
= &
- \frac {1} {2} 
\left  [ n^{-1/2} \frac { \partial \Dual (\bsm {0}, \btheta^{*}) } { \partial (\blambda, \btheta) } \right ]^{\top} 
S^{-1} 
\left [ n^{-1/2} \frac { \partial \Dual (\bsm {0}, \btheta^{*}) } { \partial (\blambda, \btheta) } \right ] 
+ o_{p} (1). 
\end{align*}  

Let  
\begin{gather*} 
\bnu_{1} = n^{-1/2} \left [ \frac { \partial \Dual (\bsm {0}, \btheta^{*}) } { \partial \blambda } \right ], \, \, \, \, \, \, 
\bnu_{2} = n^{-1/2} \left [ \frac { \partial \Dual (\bsm {0}, \btheta^{*}) } { \partial \btheta } \right ], \\ 
\Lambda = S_{\blambda \blambda} - S_{\blambda \btheta} S_{\btheta \btheta}^{-1} S_{\blambda \btheta}^{\top}, \, \, \, \, \, \, 
D = 
\begin{pmatrix} 
\mathbb {I}, & - S_{\blambda \btheta} S_{\btheta \btheta}^{-1} 
\end{pmatrix}, 
\end{gather*} 
with $ D $ and the identity matrix $ \mathbb {I} $ with proper sizes. 
We then get
\begin{align*} 
R_{n} & = 
- \bnu_{2}^{\top} S_{\btheta \btheta}^{-1} \bnu_{2} + (\bnu_{1}^{\top}, \bnu_{2}^{\top}) S^{-1} 
\begin{pmatrix}
\bnu_{1} \\ 
\bnu_{2}
\end{pmatrix} 
+ o_{p} (1) \\ 
& = \left \{ \bnu_{1} - S_{\blambda \btheta} S_{\btheta \btheta}^{-1} \bnu_{2} \right \}^{\top} 
\Lambda^{-1} 
\left \{ \bnu_{1} - S_{\blambda \btheta} S_{\btheta \btheta}^{-1} \bnu_{2} \right \} 
+ o_{p} (1) \\ 
& = 
\begin{pmatrix}
\bnu_{1} \\ 
\bnu_{2}
\end{pmatrix}^{\top} 
(D^{\top} \Lambda^{-1} D) 
\begin{pmatrix}
\bnu_{1} \\ 
\bnu_{2}
\end{pmatrix}  
+ o_{p} (1),  
\end{align*}
where the middle step can be obtained via some typical matrix algebra or Theorem 8.5.11 in \citet {harville1997matrix}. 

As given in the proof of Lemma \ref {lem:first_deriv}, the asymptotic variance of $ (\bnu_{1}, \bnu_{2}) $ is $ V $.
We also have
\begin{align*} 
D V D^{\top} 
& = 
D 
\left [
\begin{pmatrix} 
S_{\blambda \blambda} & \bsm {0} \\ 
\bsm {0} & - S_{\btheta \btheta} \\ 
\end{pmatrix} 
- S 
\begin{pmatrix} 
\bsm {0} & \bsm {0} \\ 
\bsm {0} & W \\
\end{pmatrix} 
S 
\right ]
D^{\top} \\ 
& = D 
\begin{pmatrix} 
S_{\blambda \blambda} & \bsm {0} \\ 
\bsm {0} & - S_{\btheta \btheta} \\ 
\end{pmatrix} 
D^{\top} 
- \bsm {0} \\ 
& = S_{\blambda \blambda} - S_{\blambda \btheta} S_{\btheta \btheta}^{-1} S_{\blambda \btheta}^{\top} \\ 
& = \Lambda.  
\end{align*} 
Hence,
\[ 
V (D^{\top} \Lambda^{-1} D) V (D^{\top} \Lambda^{-1} D) V = V (D^{\top} \Lambda^{-1} D) V.
\] 
By the result on quadratic forms of the multivariate normal (section 3.5, \citet {serfling2000approximation}), 
the limiting distribution of $ R_{n} $ is chi-square with the degrees of freedom being the trace of $ (D^{\top} \Lambda^{-1} D) V $, which is $ l $ as claimed in this theorem.
This completes the proof.  

\end{proof}

%%%%%%%%%%%%%%%%%%%%%%%%%%%%%%%%%%%%%%%%%%%%%%%%%%%%%%%%%%%%%%%%%%%%%%

\newpage
\section{Definability of the profile log-EL} 
\label{App:define_profile}

Discussions of the properties of the ELRT statistic are not meaningful if the profile log-EL $ \tilde \ell_{n} (\bxi) $ is not well defined. 
In fact, in some situations, the constrained maximization has no solution \citep {grendar2009empty}. 
Such an ``empty-set'' problem can be an issue, but
there are methods in the literature to overcome this obstacle \citep {chen2008adjusted, liu2010adjusted, tsao2014extended}. 
In this Appendix, we show that our $ \tilde \ell_{n} (\bxi) $ does not suffer from the ``empty-set'' problem under two additional mild conditions.
The first condition restricts our attention to quantile values $ \{ \xi_{r}: r \in I \} $ in the range  
\[
\min_{j} x_{r j} < \xi_{r} < \max_{j} x_{r j}.
\] 
The second requires one of the components of $ \bq(x) $ to be monotone in $ x $, in addition to a component being $ 1 $.
All of our examples satisfy these conditions.

To define the profile log-EL $ \tilde \ell_{n} (\bxi) $, we must have some $ p_{k j} > 0 \text { and } \btheta_{r} $ such that
\begin{gather*} 
\sum_{k, j} p_{k j} \exp (\btheta_{r}^{\top} \bq (x_{k j})) = 1, \, \, \, r = 0, 1, \ldots, m, \\ 
\sum_{k, j} p_{k j} \exp (\btheta_{r}^{\top} \bq (x_{k j})) [\mathbbm{1} (x_{k j} \leq \xi_{r}) - \tau_{r}] = 0, \, \, \, r \in I. 
\end{gather*}
We work on the most general case where $ I $ contains all populations, and without loss of generality let $ d = 2 $.
The above expressions are equivalent to (including $ r = 0 $ and allowing $ \btheta_{0} \neq 0 $)
\begin{gather*} 
\sum_{k, j} p_{k j} \exp (\btheta_{r}^{\top} \bq (x_{k j})) [\mathbbm{1} (x_{k j} \leq \xi_{r})] = \tau_{r}, \\ 
\sum_{k, j} p_{k j} \exp (\btheta_{r}^{\top} \bq (x_{k j})) [\mathbbm{1} (x_{k j} > \xi_{r})] = 1 - \tau_{r}. 
\end{gather*}  

Let
\(
\btheta_{r}^{\top} = (\theta_{r 1}, \theta_{r 2}), 
\)
and
\(
\bq^{\top} (x) = (q_{1} (x), q_{2} (x))
\)
where $ q_{1} (x) \equiv 1 $ and $ q_{2} (x) $ is monotone in $ x $. 
We can rewrite the equations as
\begin{gather*} 
\sum_{k, j} p_{k j} 
\exp \left \{ \theta_{r 1}^{'} + \theta_{r 2} [ q_{2} (x_{k j}) - q_{2} (\xi_{r}) ] \right \} [ \mathbbm{1} (x_{k j} \leq \xi_{r}) ] 
= \tau_{r}, \\ 
\sum_{k, j} p_{k j} 
\exp \left \{ \theta_{r 1}^{'} + \theta_{r 2} [ q_{2} (x_{k j}) - q_{2} (\xi_{r}) ] \right \} [ \mathbbm{1} (x_{k j} > \xi_{r}) ] 
= 1 - \tau_{r},   
\end{gather*} 
with $ \theta_{r 1}^{'} = \theta_{r 1} + \theta_{r 2} q_{2} (\xi_{r}) $. 
For notational simplicity, we retain the notation $ \theta_{r 1} $ instead of $ \theta_{r 1}^{'} $ in what follows. 

Let $ p_{k j}^{*} $ be any set of non-negative values such that 
\(
\sum_{k, j} p_{k j}^{*} = 1
\). 
Define     
\begin{align*} 
A_{r} (\theta_{r 2}) 
&=
 \sum_{k, j} p_{k j}^{*} \exp \left \{ \theta_{r 2} [q_{2} (x_{k j}) - q_{2} (\xi_{r})] \right \} 
[ \mathbbm{1} (x_{k j} \leq \xi_{r}) ] \\ 
B_{r} (\theta_{r 2}) 
& = 
\sum_{k, j} p_{k j}^{*} \exp \left \{ \theta_{r 2} [q_{2} (x_{k j}) - q_{2} (\xi_{r})] \right \} 
[ \mathbbm{1} (x_{k j} > \xi_{r}) ]. 
\end{align*}

Since $ q_{2} (x) $ is a monotone increasing function in $ x $, $ A_{r} (\theta_{r 2}) $ is decreasing in $ \theta_{r 2} $ and $ B_{r} (\theta_{r 2}) $ is increasing in $ \theta_{r 2} $. 
Thus, we have  
\begin{align*}  
\lim_{\theta_{r 2} \to - \infty} A_{r} (\theta_{r 2}) = \infty, & \, \, \, \lim_{\theta_{r 2} \to \infty} A_{r} (\theta_{r 2}) = 0; \\ 
\lim_{\theta_{r 2} \to - \infty} B_{r} (\theta_{r 2}) = 0, & \, \, \, \lim_{\theta_{r 2} \to \infty} B_{r} (\theta_{r 2}) = \infty. 
\end{align*} 
These imply that the ratio $ A_{r} (\theta_{r 2}) / B_{r} (\theta_{r 2}) $ is decreasing in $ \theta_{r 2} $ and that 
\begin{align*} 
\lim_{\theta_{r 2} \to - \infty} A_{r} (\theta_{r 2}) / B_{r} (\theta_{r 2}) = \infty, \, \, \, 
\lim_{\theta_{r 2} \to \infty} A_{r} (\theta_{r 2}) / B_{r} (\theta_{r 2}) = 0. 
\end{align*} 
By the intermediate value theorem, there must exist a value $ \theta_{r 2}^{*} $ such that 
\begin{align*} 
A_{r} (\theta_{r 2}^{*}) / B_{r} (\theta_{r 2}^{*}) = \tau_{r}/(1 - \tau_{r}). 
\end{align*} 
Let 
$
\theta_{r 1}^{*} = - \log \left \{ A_{r} (\theta_{r 2}^{*}) + B_{r} (\theta_{r 2}^{*}) \right \}
$. 
We note that $ p_{k j}^{*} $ and $ \btheta_{r}^{*} = (\theta_{r 1}^{*}, \theta_{r 2}^{*})^{\top} $ form a solution to the system. 
Hence, a solution to the system always exists. 

We may shift the solution to set $ \btheta_0 = \bsm {0} $ if required.
Validity in the general case of $ d > 2 $ is implied by setting the other entries of $ \btheta_{r} $ to the value $ 0 $.
Therefore, we have shown that our profile log-EL $ \tilde \ell_{n} (\bxi) $ does not suffer from the ``empty-set'' problem under mild conditions.

%%%%%%%%%%%%%%%%%%%%%%%%%%%%%%%%%%%%%%%%%%%%%%%%%%%%%%%%%%%%%%%%%%%%%%

\section*{Acknowledgments}
The authors are grateful to the referees and the Editor for their helpful comments and suggestions. 

\bibliographystyle{imsart-nameyear}
\bibliography{Archer_proposal}

\end{document}